\documentclass[12pt]{article}

\usepackage{amsmath,amsthm}
\usepackage{amssymb}
\usepackage{color}
\usepackage{xspace}
\usepackage{pstricks,pst-node}

\usepackage[colorlinks=true,
linkcolor=webgreen,
filecolor=webbrown,
citecolor=webgreen]{hyperref}

\definecolor{webgreen}{rgb}{0,.5,0}
\definecolor{webbrown}{rgb}{.6,0,0}

\def\red{\textcolor{red} }
\def\blue{\textcolor{blue} }
\def\v{\vert}

\def\p{\ensuremath{\mathcal P}\xspace}
\def\s{\ensuremath{\mathcal S}\xspace}

\def\o{increasing ordered tree\xspace}
\def\os{increasing ordered trees\xspace}
\def\pm{perfect matching\xspace}
\def\ps{perfect matchings\xspace}
\def\i{\ensuremath{\mathcal I}\xspace}
\def\k{\ensuremath{\mathcal K}\xspace}
\def\n{\ensuremath{\mathcal {NK}}\xspace}
\def\t{\ensuremath{\mathcal T}\xspace}
\def\gf{generating function\xspace}
\def\si{\ensuremath{\sigma}\xspace}

\def\mbf#1{\mathchoice{\hbox{\boldmath $\displaystyle #1$}}
	{\hbox{\boldmath $\textstyle #1$}}
	{\hbox{\boldmath $\scriptstyle #1$}}
	{\hbox{\boldmath $\scriptscriptstyle #1$}}} 

\newcommand{\StirlingPartition}[2]{\genfrac{ \{ }{ \} }{0pt}{}{#1}{#2}}

\begin{document}
\newtheorem{theorem}{Theorem}
\newtheorem{defn}[theorem]{Definition}
\newtheorem{lemma}[theorem]{Lemma}
\newtheorem{prop}[theorem]{Proposition}
\newtheorem{cor}[theorem]{Corollary}
\begin{center}
{\Large
Klazar Trees and Perfect Matchings                         \\ 
}

\vspace{10mm}
DAVID CALLAN  \\
Department of Statistics  \\
\vspace*{-2mm}
University of Wisconsin-Madison  \\
\vspace*{-2mm}
Medical Science Center \\
\vspace*{-2mm}
1300 University Ave  \\
\vspace*{-2mm}
Madison, WI \ 53706-1532  \\
{\bf callan@stat.wisc.edu}  \\
\vspace{4mm}

October 27, 2008
\end{center}

\begin{abstract}
Martin Klazar computed the total weight of ordered trees under 12 
different notions of weight. The 
last and perhaps most interesting of these weights, $w_{12}$, led to a recurrence relation and an identity 
for which he requested combinatorial explanations. Here we provide 
such explanations. To do so, we introduce the notion of a ``Klazar violator'' 
vertex in an increasing ordered tree and observe that $w_{12}$ 
counts what we call Klazar trees---increasing ordered trees with no Klazar violators. A highlight 
of the paper is a bijection from $n$-edge increasing ordered trees to perfect 
matchings of $[2n]=\{1,2,\ldots,2n\}$ 
that sends Klazar violators to even numbers matched to a larger odd 
number. We find the distribution of the latter matches and, in particular, establish the 
one-summation explicit formula $\sum_{k=1}^{\lfloor n/2 
\rfloor}(2k-1)!!^{2}\StirlingPartition{n+1}{2k+1}$ for the number of
perfect matchings of $[2n]$ with no even-to-larger-odd matches.
The proofs are mostly bijective.
\end{abstract}

\section{Introduction}
Martin Klazar \cite{klazar1,klazar2} defined a \emph{drawing} of an $n$-edge 
ordered tree $T$ to be a sequence of trees 
$(T_{1},T_{2},\ldots,T_{n})$ such that $T_{n}=T$ and $T_{i}$ arises 
from $T_{i-1}$ by deleting a leaf of $T_{i}$. He defined the weight 
$w_{12}$ of $T$ (the last of 12 weights he considered) to be the 
number of different drawings of $T$ and defined $w_{12}(n)=\sum_{T\in 
\t_{n}}w_{12}(T)$ where $\t_{n}$ is denotes the set of $n$-edge 
ordered trees, counted by the Catalan numbers, sequence
(\htmladdnormallink{A000108}{http://www.research.att.com:80/cgi-bin/access.cgi/as/njas/sequences/eisA.cgi?Anum=A000108})
in OEIS \cite{oeis}.
Using generating functions, he found the \gf $\sum_{n\ge 
0}w_{12}(n)x^{n}/n!=\sqrt{e^{x}/(2-e^{x})}$.
He also established the identity
\begin{equation}
    2^{n}\sum_{T\in \t_{n}} 
    w_{12}(T) \textrm{\scriptsize $\left(\frac{1}{2}\right)$}^{\ell(T)}=(2n-1)!!
    \label{eq1}
\end{equation}
where $\ell(T)$ is the number of leaves in $T$, and the recurrence relation
\begin{equation}
 w_{12}(n)=w_{12}(n-1)+\sum_{i=1}^{n-1}w_{12}(i)\binom{n-1}{i-1}, 
    \label{eq2}
\end{equation}
and wrote ``It would \ldots be interesting to give direct
combinatorial proofs and interpretations'' of (\ref{eq1}) and 
(\ref{eq2}). Here we will do so. 

The outline of the paper is as follows. 
Section 2 introduces the notion of Klazar violator and an explicit 
class of trees, Klazar trees, counted by $w_{12}$. Section 3 gives a 
combinatorial proof of identity (\ref{eq1}) and Section 4 illuminates 
it. Section 5 finds the \gf for Klazar violators. Section 6 presents 
an easily described class of perfect matchings counted by $w_{12}$, 
which serves in Section 7 to give a combinatorial interpretation of 
recurrence (\ref{eq2}). Section 8 discusses codes for trees and 
matchings that are useful as an intermediate construct for the 
bijection in the next section. Section 9 describes a bijection, 
both recursively and explicitly, between \os and perfect matchings 
represented as dot diagrams. This bijection, interesting in its own 
right, translates the result of Section 7 to provide a direct combinatorial 
proof of recurrence (\ref{eq2}) in context. Section 10 presents some 
consequences of this bijection, including a bivariate \gf for perfect matchings 
on $[2n]$ counting instances of an even number matched to a larger odd 
number and some results for the ``trapezoidal words'' considered by 
Riordan.

\section[The Weight $w_{12}$ Counts Klazar Trees]{The Weight $\protect\mbf{w_{12}}$ Counts Klazar Trees}
Let $\i_{n}$ denote the set of $n$-edge increasing ordered trees, that 
is, $n$-edge ordered trees with vertices labeled $0,1,\ldots,n$ so that the label 
of each child vertex exceeds that of its parent. It is well known, and 
indeed a nice proof is included in Klazar's paper \cite{klazar1}, 
that $\v\,\i_{n}\,\v=(2n-1)!!$, the odd double factorial. Now $w_{12}(n)$ counts trees in 
$\i_{n}$ satisfying a technical condition whose description is 
facilitated by introducing some terminology for an increasing ordered tree.
A \emph{descent} is a pair of adjacent sibling vertices in which the 
first exceeds the second. The first is a \emph{descent initiator}, the 
second a \emph{descent terminator}. 
The (left) \emph{cohort} $C(v)$ of a 
vertex $v$ is the list of all siblings of $v$ lying strictly to its left. The 
\emph{big cohort} $B(v)$ of $v$ is the maximal terminal sublist of 
the left cohort of $v$ all of whose entries exceed $v$. Thus the 
big cohort of $v$ is nonempty iff $v$ is a descent terminator. The 
\emph{associate} $A(v)$ of a descent terminator $v$ is the smallest 
entry in the big cohort of $v$. Figure 1 illustrates these notions.

\begin{center}

\begin{pspicture}(-4,-3.8)(4,3.7)
\psline(-4,2)(-2,1)(0,2)(0,3)
\psline(-2,1)(0,0)(2,1)
\psline(-3,3)(-3,2)(-2,1)(-2,2)
\psline(-1,2)(-2,1)(0,2)
\psline(-1.5,3)(-1,2)(-0.5,3)

\psdots(-3,3)(-1.5,3)(-0.5,3)(-4,2)(-3,2)(-2,2)(-1,2)(0,2)(0,3)(-2,1)(2,1)(0,0)

\rput(0,-.3){\textrm{{\footnotesize $0$}}}
\rput(-2,.7){\textrm{{\footnotesize $1$}}}
\rput(2,.7){\textrm{{\footnotesize $7$}}}

\rput(-4.2,2.2){\textrm{{\footnotesize $3$}}}
\rput(-3.2,2.2){\textrm{{\footnotesize $6$}}}
\rput(-2.2,2.2){\textrm{{\footnotesize $9$}}}
\rput(-1.3,2.2){\textrm{{\footnotesize $4$}}}
\rput(-0.2,2.2){\textrm{{\footnotesize $2$}}}

\rput(-3.2,3.2){\textrm{{\footnotesize $11$}}}
\rput(-1.7,3.2){\textrm{{\footnotesize $10$}}}
\rput(-0.6,3.2){\textrm{{\footnotesize $5$}}}
\rput(0.1,3.2){\textrm{{\footnotesize $8$}}}

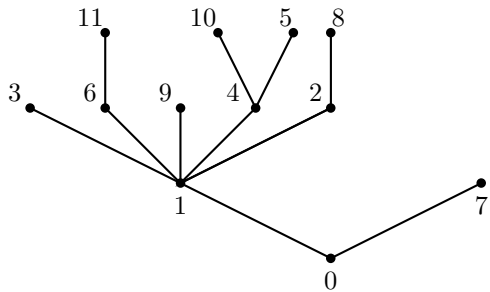
\rput(0,-1.2){\textrm{{\small An increasing ordered tree}}}  
\rput(0,-1.8){\textrm{{\small  The cohort of 1 is empty, the  
cohort of 4 is $C(4)=(3,6,9)$; 
 } }}  
\rput(0,-2.4){\textrm{{\small the big cohort of 4 is $B(4)=(6,9)$, the big cohort of 2 
is $(3,6,9,4)$; } }} 
\rput(0,-3){\textrm{{\small  the descent terminators are  4,2,5 and 
their associates are $A(4) = 6, \ A(2)= 
3,\ A(5)= 10$.  } }}

\rput(0,-3.8){\textrm{Figure 1}}

\end{pspicture}
\end{center} 

\textbf{Definition} A \emph{Klazar violator} (KV for short) in an increasing ordered 
tree is a descent terminator whose left associate 
is smaller than every child of $v$.

In particular, a descent terminator with no children is a Klazar violator 
because it satisfies this condition vacuously. For example, in the 
tree in Figure 1 above, the descent terminators are 4,2,5, and the Klazar 
violators are 2 and 5.

\textbf{Remark} If we follow the standard convention that the minimum 
of an empty set is $\infty$ and extend the notion of associate to all 
vertices, then we can say that a vertex $v$ is a  Klazar complier (KC), 
that is, not 
a Klazar violator, iff the associate of $v$ is $\ge$ the minimum of 
the children of $v$. (The weak inequality $\ge$ is used to allow for 
$\infty\ge \infty$.)

\textbf{Definition} A \emph{Klazar tree} is an increasing ordered 
tree with no Klazar violators.

An increasing (i.e. child $>$ parent) labeling of an ordered tree $T$ 
determines a drawing of $T$---delete the vertices in decreasing 
order---but the resulting drawings are not all distinct. The problem 
is that in labeling the vertices as the tree is built up from a 
drawing to produce an \o, adding an edge among a cluster of sibling 
leaf edges with a common parent gives the same tree regardless of 
where in the cluster the new edge is placed. It is straightforward to 
verify, however, that if in this situation, the new edge is always 
placed so that it is the rightmost edge of the cluster the resulting 
labeled tree will be a Klazar tree and otherwise at least one Klazar 
violator will be present. Thus $w_{12}(T)$ is the number of Klazar 
trees whose underlying ordered tree is $T$ and, letting $\k_{n}$ denote 
the set of all Klazar 
trees with $n$ edges, $w_{12}(n)=\v\,\k_{n}\,\v$.

\section{A Combinatorial Proof of Klazar's Identity}
It is convenient, following Deutsch \cite{prescribed}, to define a \emph{node} in a 
rooted tree to be a vertex that is neither the root nor a leaf. Now 
Klazar's identity (\ref{eq1}) can be written as 
\begin{equation}
    \sum_{K\in\k_{n}}2^{\nu(K)} = (2n-1)!!
    \label{eq3}
\end{equation}
where $\nu(K)$ is the number of nodes in $K$. 
Since $2^{\nu(K)}$ is 
the number of subsets of the nodes of $K$, let us define a 
\emph{node-marked} Klazar tree to be one in which some (all, or none) of its 
nodes are marked and let $\n_{n}$ denote the set of node-marked Klazar 
trees on $n$ edges. Thus (\ref{eq3}) asserts that $\v\,\n_{n}\,\v=(2n-1)!!$. To 
prove this assertion we exhibit a simple bijection $\phi$ from $\n_{n}$ to $\i_{n}$, 
the set of $n$-edge increasing ordered trees.

Given a node-marked Klazar tree $T$, turn each marked node $u$ into a 
Klazar violator by the following cut-and-paste procedure. Let $v$ be 
the smallest child of $u$. Take $v$ and its cohort and transfer all 
these vertices along with their subtrees and parent edges so that 
they become siblings of $u$ situated immediately to the left of the 
big cohort of $u$, and then remove the mark from $u$ as illustrated. 
(The $*$ superscript refers to the image tree. Thus 
$A^{*}(w),B^{*}(w),C^{*}(w)$ refer respectively to the associate, big 
cohort, and cohort of $w$ in $\phi(T)$.)
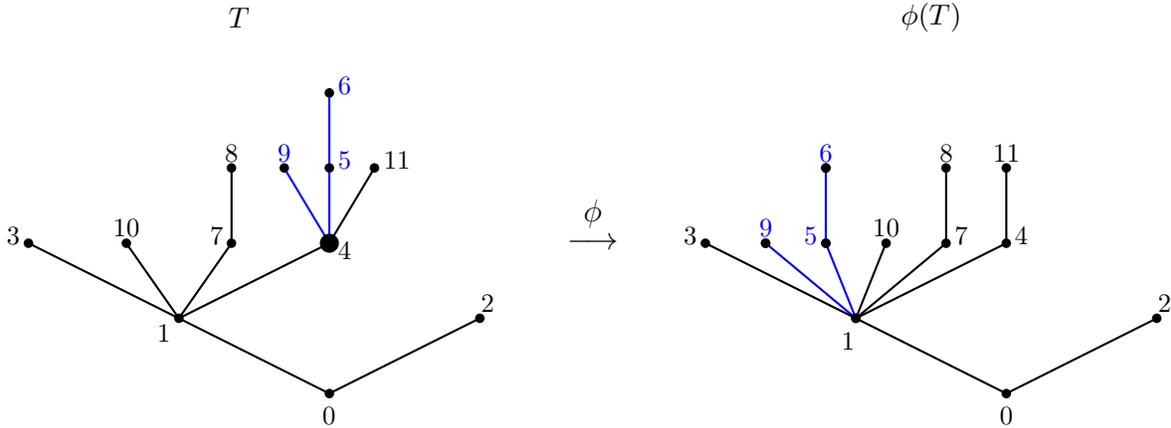
\begin{center} 
\begin{pspicture}(-4,-4.8)(10.5,5)

\psline[linecolor=blue](-0.6,3)(0,2)(0,3)(0,4)
\psline[linecolor=blue](5.8,2)(7,1)(6.6,2)(6.6,3)

\psline(-4,2)(-2,1)(0,2)(.6,3)
\psline(-2.7,2)(-2,1)(-1.3,2)(-1.3,3)
\psline(-2,1)(0,0)(2,1)

\psdots(-4,2)(.6,3)(-0.6,3)(0,2)(0,3)(0,4)(-2.7,2)(-2,1)(-1.3,2)(-1.3,3)(5.8,2)(7,1)(6.6,2)(6.6,3)(0,0)(2,1)

\rput(0,-.3){\textrm{{\footnotesize $0$}}}
\rput(-2.2,0.8){\textrm{{\footnotesize $1$}}}
\rput(2.1,1.2){\textrm{{\footnotesize $2$}}}
\rput(-4.2,2.1){\textrm{{\footnotesize $3$}}}
\rput(0.2,1.9){\textrm{{\footnotesize $4$}}}
\blue{
\rput(0.2,3.1){\textrm{{\footnotesize $5$}}}
\rput(0.2,4.1){\textrm{{\footnotesize $6$}}}
\rput(-0.6,3.2){\textrm{{\footnotesize $9$}}} }

\rput(-1.5,2.1){\textrm{{\footnotesize $7$}}}
\rput(-1.3,3.2){\textrm{{\footnotesize $8$}}}
\rput(-2.7,2.2){\textrm{{\footnotesize $10$}}}
\rput(0.9,3.1){\textrm{{\footnotesize $11$}}}

\psline(5,2)(7,1)(9,0)(11,1)
\psline(7,1)(9,2)(9,3)
\psline(7,1)(8.2,2)(8.2,3)
\psline(7,1)(7.4,2)

\psdots(5,2)(7,1)(9,0)(11,1)(5.8,2)(9,2)(9,3)(6.6,3)(6.6,2)(8.2,2)(8.2,3)(7.4,2)

\rput(3.5,2.4){$\phi$}
\rput(3.5,2){$\longrightarrow$}

\rput(9,-.3){\textrm{{\footnotesize $0$}}}
\rput(6.9,0.7){\textrm{{\footnotesize $1$}}}
\rput(11.1,1.2){\textrm{{\footnotesize $2$}}}
\rput(4.8,2.1){\textrm{{\footnotesize $3$}}}
\rput(9.2,2.1){\textrm{{\footnotesize $4$}}}
\blue{
\rput(6.4,2.1){\textrm{{\footnotesize $5$}}}
\rput(6.6,3.2){\textrm{{\footnotesize $6$}}}
\rput(5.8,2.2){\textrm{{\footnotesize $9$}}} }

\rput(8.4,2.1){\textrm{{\footnotesize $7$}}}
\rput(8.2,3.2){\textrm{{\footnotesize $8$}}}
\rput(7.4,2.2){\textrm{{\footnotesize $10$}}}
\rput(9,3.2){\textrm{{\footnotesize $11$}}}

\rput(-1.2,5){\textrm{{\small $T$}}} 
\rput(8,5){\textrm{{\small $\phi(T)$}}}

\rput(0,-1.6){\textrm{{\small $u=4$ is a marked node in $T$,}}} 
\rput(0,-2.2){\textrm{{\small $v=5$ is smallest child of $u$,}}} 
\rput(0,-2.8){\textrm{{\small $(9)$ is the cohort of $v$.}}} 

\rput(8,-1.2){\textrm{{\small $u=4$ is a Klazar violator in $\phi(T)$,}}} 
\rput(8,-2.2){\textrm{{\small $B^{*}(u)=
(\:\underbrace{9\vphantom{/}}_{\textrm{{\footnotesize $C(v)$}}},\,
\underbrace{5\vphantom{|}}_{\textrm{{\footnotesize 
$\vphantom{|}v$}}},\,\underbrace{10,\,7}_{\textrm{{\footnotesize $B(u)$}}}\:),$}}} 
\rput(8,-3.2){\textrm{{\small $B^{*}(v)=C(v)=(9).$}}} 
\rput(3.5,-4.2){\textrm{{\small  In $\phi(T)$, the marked node $u$ is recaptured  as a Klazar violator,
the vertex $v$ is recaptured }}} 
\rput(3.5,-4.8){\textrm{{\small  as $A^{*}(u)$ because $v$ is the 
smallest entry in $B^{*}(u)$, and $C(v)$ is recaptured 
as $B^{*}(v)$.   }}}

\psset{dotscale=2} \psdots(0,2)

\end{pspicture}
\end{center} 
In this example there is just one marked node; if there is more than 
one, the marked nodes are processed one after the other (the order of 
processing is immaterial) as in the following example.

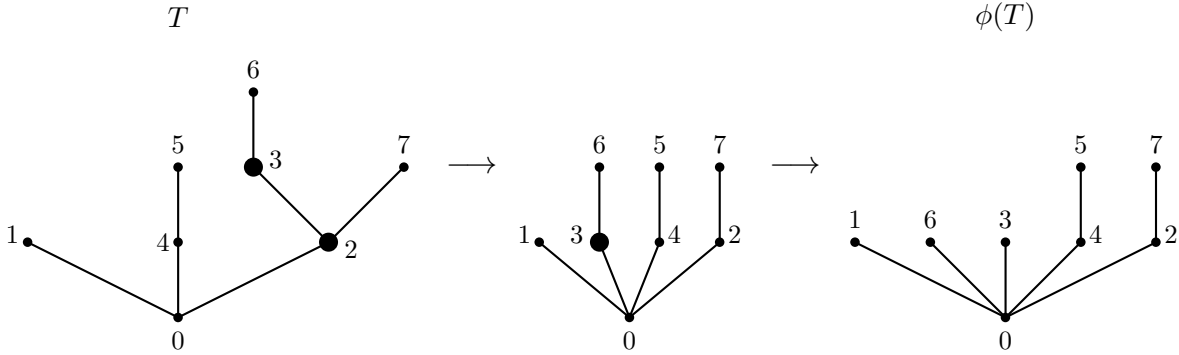
\begin{center} 
\begin{pspicture}(-8,-.5)(7,4.3)

\psline(-8,1)(-6,0)(-6,1)(-6,2)
\psline(-6,0)(-4,1)(-5,2)(-5,3)
\psline(-4,1)(-3,2)

\psdots(-8,1)(-6,0)(-6,1)(-6,2)(-4,1)(-5,2)(-5,3)(-4,1)(-3,2)

\rput(-6,-.3){\textrm{{\footnotesize $0$}}}
\rput(-8.2,1.1){\textrm{{\footnotesize $1$}}}
\rput(-3.7,.9){\textrm{{\footnotesize $2$}}}
\rput(-4.7,2.1){\textrm{{\footnotesize $3$}}}
\rput(-6.2,1){\textrm{{\footnotesize $4$}}}
\rput(-6,2.3){\textrm{{\footnotesize $5$}}}
\rput(-5,3.3){\textrm{{\footnotesize $6$}}} 
\rput(-3,2.3){\textrm{{\footnotesize $7$}}}

\psline(-1.2,1)(0,0)(1.2,1)(1.2,2)
\psline(-.4,2)(-.4,1)(0,0)(.4,1)(.4,2)

\psdots(-1.2,1)(0,0)(1.2,1)(1.2,2)(-.4,2)(-.4,1)(.4,1)(.4,2)

\rput(-2.1,2){$\longrightarrow$}

\rput(0,-.3){\textrm{{\footnotesize $0$}}}
\rput(-1.4,1.1){\textrm{{\footnotesize $1$}}}
\rput(1.4,1.1){\textrm{{\footnotesize $2$}}}
\rput(-.7,1.1){\textrm{{\footnotesize $3$}}}
\rput(.6,1.1){\textrm{{\footnotesize $4$}}}
\rput(.4,2.3){\textrm{{\footnotesize $5$}}}
\rput(-.4,2.3){\textrm{{\footnotesize $6$}}} 
\rput(1.2,2.3){\textrm{{\footnotesize $7$}}}

\rput(-6,4){\textrm{{\small $T$}}} 
\rput(5,4){\textrm{{\small $\phi(T)$}}}

\rput(2.2,2){$\longrightarrow$}

\psline(3,1)(5,0)(7,1)(7,2)
\psline(4,1)(5,0)(6,1)(6,2)
\psline(5,0)(5,1)

\psdots(3,1)(5,0)(7,1)(7,2)(4,1)(6,1)(6,2)(5,1)

\rput(5,-.3){\textrm{{\footnotesize $0$}}}
\rput(3,1.3){\textrm{{\footnotesize $1$}}}
\rput(7.2,1.1){\textrm{{\footnotesize $2$}}}
\rput(5,1.3){\textrm{{\footnotesize $3$}}}
\rput(6.2,1.1){\textrm{{\footnotesize $4$}}}
\rput(6,2.3){\textrm{{\footnotesize $5$}}}
\rput(4,1.3){\textrm{{\footnotesize $6$}}} 
\rput(7,2.3){\textrm{{\footnotesize $7$}}}

\psset{dotscale=2} \psdots(-5,2)(-4,1)(-.4,1)

\end{pspicture}
\end{center} 

As each marked node is processed, it is turned into a Klazar violator 
without affecting the complier/violator status of any other vertex. 
Thus the originally marked nodes can be recovered as the Klazar 
violators in $\phi(T)$, and the entire process is reversible.

\section{An Interpretation of Klazar's Identity}

From the preceding section we have a combinatorial proof of (\ref{eq1}) but not yet 
a satisfactory combinatorial interpretation: what is still missing is 
a nice characterization of the image in $\i_{n}$, under $\phi$, of the 
node-marked Klazar trees with $\ell$ leaves. Such a characterization 
is far from obvious, so how to find one?
We will tackle this problem (with gratifying success) but first let 
us review the genesis of (\ref{eq1}). Klazar, looking for the \gf for 
$w_{12}(n)$, found a recurrence for a more refined count by number of 
leaves. This led him to to an expression for the bivariate \gf 
\begin{equation}
    F^{*}(x,y) := \sum_{n\ge 0}\sum_{T\in \t_{n}}w_{12}(T) 
    \frac{x^{n}}{n!}y^{\ell(T)},
    \label{fstar}
\end{equation}
where $\ell(T)$ is the number of leaves of $T$, to wit,
\begin{equation}
 F^{*}(x,y) = \sqrt{ \frac{2y-1}{2ye^{x(1-2y)}-1} }.
 \label{leaves}
\end{equation}
Of course, $F^{*}(x,1)=\sqrt{e^{x}/(2-e^{x})}$ is the desired exponential 
\gf (egf) 
for $w_{12}(n)$. But Klazar also noted that the tweaked function
$F^{**}(x,y):=F^{*}(2x,y/2)$ has the curious property that 
$F^{**}(x,1)=1/\sqrt{ 1-2x }$, the exponential \gf for the odd double 
factorials. Since (\ref{fstar}) says
\begin{equation}
F^{**}(x,y)=\sum_{n\ge 0}2^{n}\sum_{T\in \t_{n}}w_{12}(T) 
\textrm{\footnotesize $\left(\frac{1}{2}\right)$}^{\ell(T)}
\frac{x^{n}}{n!}y^{\ell(T)},
 \label{fstar1}
\end{equation}
Klazar obtained (\ref{eq1}) by setting $y=1$ in (\ref{fstar1}) and 
equating coefficients of $\frac{x^{n}}{n!}$. 

Now let $a(n,\ell)$ denote the coefficient of $\frac{x^{n}}{n!}y^{\ell}$ in 
$F^{**}(x,y) = \sqrt{ \frac{\textrm{ \raisebox{0.5ex}{$y-1$} }}{y e^{2x(1-y)}-1} } = 1 + \sum_{n,\ell\ge 
1}a(n,\ell)\frac{x^{n}}{n!}y^{\ell}$. 
The first few values of $a(n,\ell)$ are given in the following table.

\vspace*{-4mm}

\[
\begin{array}{c|ccccccc}
	n^{\textstyle{\,\backslash \,\ell}}  & 1 & 2 & 3 & 4 & 5 & 6 & 7 \\
\hline 
	1&    1 &   & & &  &  & \\
 	2&    2 & 1 & & &  &  &   \\
	3&    4 & 10  &1 & &  &  &  \\
	4&    8 & 60  & 36 &1 &   &   & \\ 
	5&    16 & 296 & 516 & 116 &  1  &  & \\
	6&    32 & 128 & 5158 & 3508 &  358  & 1 &  \\
	7&    64 & 5664 & 42960 & 64240 & 21120 & 1086 & 1 
\end{array}
\]

\centerline{\textrm{Table of values of $a_{n,\ell}$}}

The first order of business is to try to find a statistic 
on $\i_{n}$ whose distribution is given by the array 
$\big(a(n,\ell)\big)$. The first column 
$(a(n,1))_{n\ge 1}$ appears to be $(2^{n-1})_{n\ge 1}$ and $2^{n-1}$ is the 
number of compositions of $n$. A composition $n=n_{1}+n_{2}+\ldots 
+n_{r}\ (n_{i}\ge 1)$ suggests an increasing ordered tree in a simple 
way: split $[n]$ into blocks of consecutive integers of lengths 
$n_{1},n_{2},\ldots,n_{r}$, and use the blocks as sibling lists, each 
having as common parent the last entry of the previous list (or the 
root, in the case of the first list) as 
illustrated in Figure 2. 
\begin{center}

\begin{pspicture}(-2,-2)(4,4)
\psline(-2,1)(0,0)(2,1)
\psline(0,0)(0,1)

\psline(1,2)(2,1)(3,2)
\psline(2,3)(3,2)(3,3)
\psline(3,2)(4,3)(4,4)

\psdots(-2,1)(0,0)(0,1)(2,1)(1,2)(3,2)(2,3)(3,3)(4,3)(4,4)

\rput(0,-.3){\textrm{{\footnotesize $0$}}}
\rput(-2,1.2){\textrm{{\footnotesize $1$}}}
\rput(0,1.2){\textrm{{\footnotesize $2$}}}

\rput(2.3,1){\textrm{{\footnotesize $3$}}}
\rput(1,2.2){\textrm{{\footnotesize $4$}}}
\rput(3.3,2){\textrm{{\footnotesize $5$}}}
\rput(2,3.2){\textrm{{\footnotesize $6$}}}
\rput(3,3.2){\textrm{{\footnotesize $7$}}}
\rput(4.2,3.2){\textrm{{\footnotesize $8$}}}
\rput(4,4.3){\textrm{{\footnotesize $9$}}}

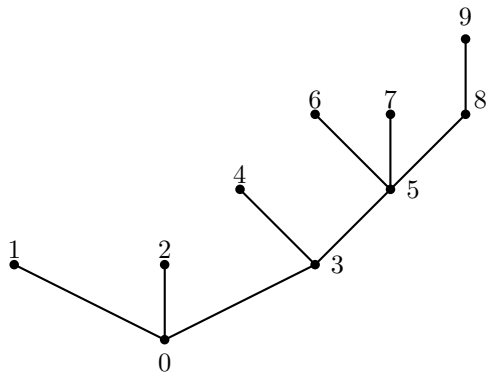
\rput(0,-1.2){\textrm{{\small The increasing ordered tree 
corresponding to the composition $(3,2,3,1)$ of 9}}}

\rput(0,-2){\textrm{Figure 2}}

\end{pspicture}
\end{center} 
These trees, counted by $2^{n-1}$, are clearly increasing, and they 
are characterized by the further properties:

\vspace*{-5mm}

\begin{itemize}
    \item  no sibling descents 
\vspace*{-2mm}
    \item  only the rightmost child of a vertex can have children.
\end{itemize}

\vspace*{-5mm}

\noindent This motivates us to define a \emph{bad} vertex in an increasing ordered 
tree to be a vertex that (i) initiates a sibling descent or (ii) initiates a sibling ascent 
and has children. Thus the $n$-edge increasing ordered trees counted by 
$a(n,1)=2^{n-1}$ are those with no bad vertices. Could it be that $a(n,\ell)$ is 
the number with $\ell-1$ bad vertices? 
Computer calculations suggest that indeed it is, and so we are 
(strongly) motivated to check if $\phi$ sends 
node-marked Klazar trees with $\ell$ leaves to increasing ordered trees with $\ell-1$ bad 
vertices. It doesn't, but it does send them to increasing ordered trees with $\ell-1$ \emph{reverse-bad} 
vertices (a vertex is 
reverse-bad if it is bad viewing the tree from right to left, that is, if 
it is bad in the tree obtained by flipping the original tree over a 
vertical line). This follows from the following two key observations.
\begin{prop}
    $($i$\,)$ In a Klazar tree, 
the number of reverse-bad vertices is 
one less than the number of leaves, and
$($ii$\,)$ the bijection $\phi:\n_{n}  \rightarrow \i_{n}$ presented above preserves the number of reverse-bad 
vertices.
\end{prop}
\textbf{Proof}\quad (i) Let $T$ be a Klazar tree. 
Given a leaf $u$ in $T$, consider the path from $u$ to the root. Let 
$\pi(u)$ be the first vertex on this path (possibly $u$ itself) that 
has a left sibling. The map $\pi$ is defined for all leaves except the 
that terminates the leftmost path from the root. We claim it is a 
bijection to the reverse-bad vertices of $T$ and the result follows. 
To see the claim, observe that if $\pi(u)=u$ then $u$ has closest 
left sibling $v$ and $v<u$ for otherwise $u$, being a leaf, would be a 
Klazar violator. Hence $u$ is reverse-bad. If $\pi(u)\ne u$ 
then $\pi(u)$ has both a left sibling and a child and so is certainly 
reverse-bad. Thus $\pi$ sends all but the exceptional leaf to 
reverse-bad vertices. Conversely, given a reverse-bad vertex $v$, map 
it to the leaf terminating the leftmost path from $v$ away from the 
root. This map is the inverse of $\pi$.

(ii) It suffices to verify the assertion for a single application of 
the ``mark to violator'' process and this involves a routine check of 
various cases, which we leave to the 
reader. \qed

Of course, the distribution of 
reverse-bad vertices is the same as the distribution of bad vertices.
The preceding discussion shows how I both stumbled upon and proved the following combinatorial 
interpretation for the array $\big(a(n,\ell)\big)$ implicit in (\ref{eq1}).
\begin{theorem}
    $a(n,\ell)$ is the number of $n$-edge increasing ordered trees with $\ell$ bad vertices 
    where a vertex is bad iff it either initiates a sibling descent or initiates a sibling ascent 
    and has children. \qed
\end{theorem}

\section{The \gf for Klazar violators}

The following result sheds further light on (\ref{eq1}).
\begin{theorem}
    Let $F(x,y,z)$ denote the trivariate \gf $\sum_{n,i,j\ge 
    0}a_{n,i,j}\frac{x^{n}}{n!}y^{i}z^{j}$ where $a_{n,i,j}$ is the number of 
    \os on $[n]$ with $i$ Klazar violators and $j$ leaves that do not 
    terminate a sibling descent. Then
    \[
    F(x,y,z) =  \left( \frac{1 + y - 2z}{1 + y - 2ze^{x(1 + y - 
    2z})}\right)^{\frac{1}{2}}.
    \label{trivariate}
    \]
\end{theorem}
We defer the proof to list some simple corollaries.
\begin{cor}
    \label{kvgf}
    The \gf for \os by number of Klazar violators is
    \[
     \left( \frac{1 - y}{2e^{x(y - 1)} - 1 - y } \right)^{\frac{1}{2}}.
     \]
\end{cor}
Proof. Put $z=1$ in $F(x,y,z)$.
\begin{cor}
    The statistics ``\#\: 
    non-descent-terminator leaves'' and ``\#\: 
    reverse-bad vertices'' are equidistributed on \os.
\end{cor}
Proof. The \gf for the first of these statistics is $F(x,1,z)=
\sqrt{ \frac{\textrm{ \raisebox{0.3ex}{$z-1$} }}{z e^{2x(1-z)}-1} }$ 
and this agrees with $F^{**}$ above.

We can also recover Klazar's bivariate \gf (\ref{leaves}) by setting  
$y=0$ in $F(x,y,z)$: \emph{every} leaf in a Klazar tree is 
a non-descent-terminator leaf since a descent-terminator leaf would be a 
Klazar violator.

\textbf{Proof of Theorem \ref{trivariate}}\quad Consider the effect of 
adding a leaf $n$ to an \o of size $n-1$. On the one hand, the number of Klazar 
violators increases by 1 if the new leaf is the immediate left 
sibling of what was originally a non-descent-terminator leaf; 
otherwise it stays the same. On the other hand, the number of non-descent-terminator leaves 
stays the same if the new leaf is either the rightmost child or the immediate left 
sibling of what was originally a non-descent-terminator leaf; 
otherwise it increases by 1. These observations lead to the recurrence 
relation
\[
a_{n,i,j}=ja_{n-1,i,j}+ja_{n-1,i-1,j}+(2n-2j+1)a_{n-1,i,j-1}
\]
for $n\ge 1,\ i\ge 0,\ j\ge 1$ and $(n,i,j)$ neither $(1,0,2)$ nor $(1,1,1)$,
with initial conditions
$a_{0,0,1}=1,\ a_{1,0,2} = 0,\ a_{1,1,1} = 0,\
a_{0,i,j}=0$ for $(i,j)\ne (0,1)$ and $a_{n,i,j}=0$ if $i<0$ or $j<1$.

This recurrence translates to the first-order partial differential 
equation
\[
(2x z - 1)F_{x}+(z + y z - 2 z^2)F_{z} + z F=0
\]
with solution (\`{a} la \cite[p.\,207]{klazar1}) as asserted in the Theorem. \qed

By similar considerations it is also possible to obtain a recurrence relation for 
$a_{n,i,j,k}$, the number of 
\os on $[n]$ with $i$ Klazar violators, $j$ leaves that do not 
terminate a sibling descent, and $k$ leaves altogether:
\[
\begin{array}{c}
    a_{0,0,1,1}=a_{1,0,1,1}=a_{2,0,2,2}=a_{2,0,1,1}=a_{2,1,1,2}=1  \textrm{ 
    and for other } n\le 2,\ a_{n,i,j,k} = 0,  \\[2mm]
    \textrm{and for } n\ge 3,\ a_{n,i,j,k} =  \hspace*{100mm} \\
     j a_{n-1,i,j,k} + j a_{n-1,i-1,j,k-1} + 
   (k-(j-1))a_{n-1,i,j-1,k} + (2n-k-j+1)a_{n-1,i,j-1,k-1}.
\end{array}
\]

The presence of the statistic ``total number of leaves'', however, precludes 
finding an ``elementary'' \gf because the known (marginal) distribution 
of this statistic is not elementary; see comment on sequence
\htmladdnormallink{A008517}{http://www.research.att.com:80/cgi-bin/access.cgi/as/njas/sequences/eisA.cgi?Anum=A008517}
in \cite{oeis}.

\section{A Class of Perfect Matchings}
A \pm (always on the support set $[2n]=\{1,2,\ldots,2n\}$) is a partition of $[2n]$ into 
2-element subsets or \emph{matches}. The \emph{size} of the matching is $n$ 
and we write all matches with the 
smaller entry first so that, for example, an even-to-odd match is an 
instance of an even number matched to a larger odd number. Thus the \pm 
$1\,5\,/\,2\,7\, / 3\,4\,/\,6\,8$ has one 
even-to-odd match, namely 2\,7. We will show that the number $a(n)$ 
of \ps of size $n$ with no even-to-odd matches 
has the same \gf as $w_{12}(n)$---$\sum_{n\ge 
0}a(n)x^{n}/n!=\sqrt{e^{x}/(2-x)}$---by finding an explicit formula 
for $a(n)$.
\begin{prop}
The number of  perfect matchings on $\{1,2,\ldots,2n\}$
with no even-to-odd matches and $k$ even-to-even matches is 
$(2k-1)!!^{2}\StirlingPartition{n+1}{2k+1}$.
\label{pmcnt}
\end{prop}
Here $\StirlingPartition{n}{k}$ is the Stirling partition number: the number 
of partitions of $[n]:=\{1,2,\ldots,n\}$ into $k$ nonempty disjoint 
sets (blocks). We defer the proof of Prop. \ref{pmcnt} to deduce the 
\gf. By convention, $(-1)!!=1$.
\begin{prop}
    \[
    \sum_{n\ge 0}\left(\sum_{k \ge 0}(2k-1)!!^{2} 
    \StirlingPartition{n+1}{2k+1} \right) \frac{x^{n}}{n!} = 
    \sqrt{\frac{e^{x}}{2-e^{x}}}.
    \]
\end{prop}
\noindent \textbf{Proof}\quad Reversing the order of summation, the left hand 
side is
\begin{align}
   & D_{x}\left(   \sum_{k \ge 0}(2k-1)!!^{2} \sum_{n\ge 0}
    \StirlingPartition{n+1}{2k+1}  \frac{x^{n+1}}{(n+1)!} \right)  
    \notag \\
   & =  D_{x}\left(   \sum_{k \ge 0}(2k-1)!!^{2} \sum_{n\ge 1}
    \StirlingPartition{n}{2k+1}  \frac{x^{n}}{n!} \right) \notag \\
  &  =  D_{x}\left(   \sum_{k \ge 0}(2k-1)!!^{2} 
    \frac{(e^{x}-1)^{2k+1}}{(2k+1)!} \right)  \notag \\
   & = e^{x} \sum_{k\ge 0} (2k-1)!!^{2} \frac{(e^{x}-1)^{2k}}{(2k)!} 
   \notag \\
    & = e^{x} \sum_{k\ge 0} \binom{2k}{k} \left( \frac{e^{x}-1}{2} 
    \right)^{2k}  \notag \\
    & = e^{x} \left( 1-(e^{x}-1)^{2} \right)^{-1/2} \notag \\
  & = \sqrt{\frac{e^{x}}{2-e^{x}}},  \notag
\end{align}
where a standard generating function for Stirling numbers \cite[Eq. 
(7.49), p.\,351]{gkp} is used at the 
second equality. \qed

 \textbf{Remark}\quad A similar calculation gives a bivariate 
generating function. If 
there are $k$ even-to-even matches, then there are also $k$ odd-to-odd matches 
and so $n-2k$ opposite-parity matches. If $a(n,j)$ denotes the number of no-even-to-odd  
matchings with $j$ odd-to-even matches, then the mixed generating 
function for $a(n,j)$ is given by
\[
\sum_{n,j\ge 0}a(n,j)\frac{x^{n}}{n!}y^{j}=
\frac{ye^{xy}}{\sqrt{ y^{2}-1+e^{xy}(2-e^{xy})}   } 
\]
\qed

It will be convenient to represent a \pm as a dot diagram with 
vertices arranged in two rows as illustrated.
\begin{center} 
    
\begin{pspicture}(-5,-1)(0,1.5)


\psline(-3,0)(-3,1)
\psline(-2,0)(0,1)
\psline(-1,1)(0,0)

\psdots(-4,1)(-3,1)(-2,1)(-1,1)(0,1)
\psdots(-4,0)(-3,0)(-2,0)(-1,0)(0,0)

\rput(-4,1.2){\textrm{{\scriptsize $1$}}}
\rput(-3,1.2){\textrm{{\scriptsize $3$}}}
\rput(-2,1.2){\textrm{{\scriptsize $5$}}}
\rput(-1,1.2){\textrm{{\scriptsize $7$}}}
\rput(0,1.2){\textrm{{\scriptsize $9$}}}

\rput(-4,-0.2){\textrm{{\scriptsize $2$}}}
\rput(-3,-0.2){\textrm{{\scriptsize $4$}}}
\rput(-2,-0.2){\textrm{{\scriptsize $6$}}}
\rput(-1,-0.2){\textrm{{\scriptsize $8$}}}
\rput(0,-0.2){\textrm{{\scriptsize $10$}}}

\psbezier[linewidth=.8pt](-4,1)(-3.4,1.5)(-2.6,1.5)(-2,1)
\psbezier[linewidth=.8pt](-4,0)(-3,0.5)(-2,0.5)(-1,0)

\rput(-2.2,-0.9){\textrm{{\small a perfect-matching (PM) dot diagram}}}

\end{pspicture}
\end{center} 
It is also convenient to distinguish an \emph{arc} joining two dots in 
the same row (same parity matches) and a \emph{line} joining dots in 
different rows (opposite parity matches). An even-to-odd match shows 
up in the dot diagram as an \emph{upline} (line of positive slope) 
and an odd-to-even match as a \emph{weak downline} (a vertical line 
or line of negative slope). Thus 6-9 is an upline, 1-5 is an arc, 
and 3-4 is a weak downline. The labels are not necessary in a PM dot 
diagram and we will often use $i$ bot to refer to the $i$th dot in the 
bottom row and analogously for $i$ top.

To establish Prop. \ref{pmcnt} bijectively, we will actually give a more 
refined count of the no-even-to-odd matchings and then invoke the identity
\[
\StirlingPartition{n+1}{2k+1} =\sum_{j\ge 0} 
\StirlingPartition{j}{2k}(2k+1)^{n-j},
\]
which is easily proved bijectively \cite[p.\,106, Identity 
201]{art}: the right hand side counts the partitions of $[n+1]$ into 
$2k+1$ blocks by smallest entry of the last block where, 
as throughout this note, the blocks of a partition 
are arranged in a standard order so that the smallest entries are increasing left 
to right.

\begin{theorem}
    \label{main}
    The number of perfect matchings on the set $[2n]$ in which no 
    even number is matched to a larger odd number, with $k$ even-to-even 
    matches and $2j$ the largest number occurring among the even-to-even 
    matches is $ (2k-1)!!^{2}\StirlingPartition{j}{2k}(2k+1)^{n-j}$.
\end{theorem}

Before proceeding with the proof we establish simple combinatorial 
interpretations of the Stirling partition numbers and the perfect 
powers in terms of matchings in dot diagrams.

\textbf{Definition}\quad An $(n,k)$ \emph{Stirling matching} is a $2\times n$ array of 
    dots with $n-k$ disjoint edges, each connecting a dot in the top 
    row to a dot lying strictly to its right in the bottom row.

\begin{center} 
\begin{pspicture}(-4,-2.5)(2,1.2)
\psset{xunit=1.2cm,yunit=.8cm}

\psline(-4,1)(0,-0.5)
\psline(-2,1)(-1,-0.5)
\psline(-1,1)(1,-0.5)
\psdots(-4,1)(-3,1)(-2,1)(-1,1)(0,1)(1,1)(2,1)
\psdots(-4,-0.5)(-3,-0.5)(-2,-0.5)(-1,-0.5)(0,-0.5)(1,-0.5)(2,-0.5)
{\gray 
\rput(-4,-0.8){\textrm{{\scriptsize $1$}}}
\rput(-3,-0.8){\textrm{{\scriptsize $2$}}}
\rput(-2,-0.8){\textrm{{\scriptsize $3$}}}
\rput(-1,-0.8){\textrm{{\scriptsize $4$}}}
\rput(0,-0.8){\textrm{{\scriptsize $5$}}}
\rput(1,-0.8){\textrm{{\scriptsize $6$}}}
\rput(2,-0.8){\textrm{{\scriptsize $7$}}} }
 
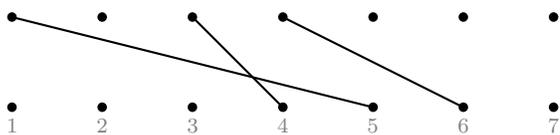
\rput(-1,-1.5){\textrm{{\footnotesize 7 dots in each row, 4 unmatched dots in 
each row}}}

\rput(-1,-2.5){\textrm{{\small A (7,4) Stirling matching}}}

\end{pspicture}
\end{center} 

\begin{prop}
    The number of $(n,k)$ Stirling matchings is the Stirling partition number
    $\StirlingPartition{n}{k}$. \label{p1}
\end{prop}
\noindent \textbf{Proof}\quad Let $\s(n,k)$ denote the set of $(n,k)$ Stirling matchings 
and set $S(n,k)=\v\,\s(n,k)\,\v$. A matching in $\s(n,k)$ can be 
obtained either (i) from a matching in $\s(n-1,k-1)$ by adding an 
unmatched dot at the end of each row---$S(n-1,k-1)$ choices---or (ii) 
from one in $\s(n-1,k)$ by similarly adding two dots and connecting 
the bottom one to any of the $k$ unmatched dots in the original top 
row---$k S(n-1,k)$ choices. Every matching in $\s(n,k)$ arises 
uniquely in one of these two ways, and so $S(n,k)$ satisfies the 
basic recurrence for the Stirling partition numbers: $S(n,k) 
=S(n-1,k-1)+k S(n-1,k)$. Since the boundary conditions 
$S(1,1)=\StirlingPartition{1}{1}=1$ and 
$S(n,k)=\StirlingPartition{n}{k}=0$ for $k>n$ or $k<1$ also hold, we 
conclude that $S(n,k)=\StirlingPartition{n}{k}$. \qed

\textbf{Remark}\quad It is possible to derive from this recurrence a 
bijection from $\s(n,k)$ to the partitions of $[n]$ into $k$ blocks 
(arranged in standard order: smallest entries increasing left 
to right). First label the dots in each row $1,2,\ldots,n$ left to 
right. Then for each $i\in [n]$, place $i$ in block $j$ ($1 \le j 
\le k$) as follows. If dot $i$ in the bottom row is unmatched, then 
$j=1+$ \#\,unmatched dots in the bottom row with label $<i$. If dot 
$i$ in the bottom row is matched, say to dot $k$ in the top row, then 
$j=k-$ \#\,edges that connect a dot $<k$ in the top row to a dot $<i$ 
in the bottom row. For example, the Stirling partition illustrated 
above corresponds to the partition
1\,5/2\,6/3\,4/7. \qed

\textbf{Definition}\quad A $(k,n)$ perfect-power matching is a 2-row 
array consisting of $k+n$ dots in the top row and $n$ dots flush right 
in the bottom row and a matching of all the lower dots to upper dots 
such that the bottom dot of each edge lies strictly to the 
right of its top dot.

\begin{center} 
\begin{pspicture}(-4,-1.8)(3,1.2)
\psset{xunit=1.2cm,yunit=.8cm}

\psline(-3,1)(2,-0.5)
\psline(-1,1)(1,-0.5)
\psline(2,1)(3,-0.5)
\psdots(-4,1)(-3,1)(-2,1)(-1,1)(0,1)(1,1)(2,1)(3,1)
\psdots(1,-0.5)(2,-0.5)(3,-0.5)
{\gray 
\rput(-4,1.3){\textrm{{\scriptsize $1$}}}
\rput(-3,1.3){\textrm{{\scriptsize $2$}}}
\rput(-2,1.3){\textrm{{\scriptsize $3$}}}
\rput(-1,1.3){\textrm{{\scriptsize $4$}}}
\rput(0,1.3){\textrm{{\scriptsize $5$}}}
\rput(1,-0.8){\textrm{{\scriptsize $1$}}}
\rput(2,-0.8){\textrm{{\scriptsize $2$}}}
\rput(3,-0.8){\textrm{{\scriptsize $3$}}}}

\rput(-0.5,-1.8){\textrm{{\small A (5,3) perfect-power matching}}}

\end{pspicture}
\end{center} 
\begin{prop}
    The number of $(k,n)$ perfect-power matchings is $k^{n}$. 
    \label{p2}
\end{prop}
\noindent \textbf{Proof}\quad Every $(k,n)$ matching comes from a $(k,n-1)$ 
matching by appending a dot to each row and connecting the lower 
dot to one of the $k$ unmatched dots in the original top row. This gives a 
multiplying factor of $k$ each time $n$ is incremented, and the result 
follows. \qed

\noindent \textbf{Proof of Theorem \ref{main}}\quad Let us take, as a working 
example, the matching
\[
1 \ 2 \:/\: 3 \ 15\:/\: 4 \ 8 \:/\: 5 \ 14 \:/\: 6 \ 12 \:/\: 7 \ 10 \:/\: 9 \ 13 \:/\: 11 \ 16
\]
with $n=8,\ k=2$ and $j=6$.

First, we take care of the $(2k-1)!!^{2}$ factor (this step is easy). 
Pick out the $k$ pairs consisting of two even integers, here 4\:8 and 
6\:12 and let $A$ denote their support, here $\{4,6,8,12\}$. These 
pairs form a perfect matching on $A$---$(2k-1)!!$ possibilities---and 
so we may extract a $(2k-1)!!$ factor and assume the pairs in question 
form a standard matching on $A$: smallest entry of $A$ matched to next 
smallest, third smallest to fourth smallest and so on. Likewise for 
the odd-to--odd pairs (also necessarily $k$ in number) we may extract 
another $(2k-1)!!$ factor and assume 
a standard matching on their support, $B$. Standardizing the matchings 
on $A$ and $B$ for our working example, we get a PM dot diagram:

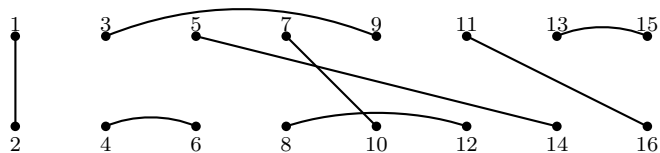
\begin{center} 
    
\begin{pspicture}(-4,-0.7)(3,1.6)
\psset{xunit=1.2cm,yunit=.8cm}

\psline(-4,1)(-4,-0.5)
\psline(-2,1)(2,-0.5)
\psline(-1,1)(0,-0.5)
\psline(1,1)(3,-0.5)
\psdots(-4,1)(-3,1)(-2,1)(-1,1)(0,1)(1,1)(2,1)(3,1)
\psdots(-4,-0.5)(-3,-0.5)(-2,-0.5)(-1,-0.5)(0,-0.5)(1,-0.5)(2,-0.5)(3,-0.5)
 
\rput(-4,1.2){\textrm{{\scriptsize $1$}}}
\rput(-3,1.2){\textrm{{\scriptsize $3$}}}
\rput(-2,1.2){\textrm{{\scriptsize $5$}}}
\rput(-1,1.2){\textrm{{\scriptsize $7$}}}
\rput(0,1.2){\textrm{{\scriptsize $9$}}}
\rput(1,1.2){\textrm{{\scriptsize $11$}}}
\rput(2,1.2){\textrm{{\scriptsize $13$}}}
\rput(3,1.2){\textrm{{\scriptsize $15$}}}

\rput(-4,-0.8){\textrm{{\scriptsize $2$}}}
\rput(-3,-0.8){\textrm{{\scriptsize $4$}}}
\rput(-2,-0.8){\textrm{{\scriptsize $6$}}}
\rput(-1,-0.8){\textrm{{\scriptsize $8$}}}
\rput(0,-0.8){\textrm{{\scriptsize $10$}}}
\rput(1,-0.8){\textrm{{\scriptsize $12$}}}
\rput(2,-0.8){\textrm{{\scriptsize $14$}}}
\rput(3,-0.8){\textrm{{\scriptsize $16$}}}

\psbezier[linewidth=.8pt](-3,1)(-2,1.6)(-1,1.6)(0,1)
\psbezier[linewidth=.8pt](2,1)(2.3,1.2)(2.7,1.2)(3,1)
\psbezier[linewidth=.8pt](-3,-0.5)(-2.7,-.3)(-2.3,-.3)(-2,-0.5)
\psbezier[linewidth=.8pt](-1,-0.5)(-0.4,-.2)(0.4,-.2)(1,-0.5)

\end{pspicture}
\end{center} 

The defining characteristics of the PM dot diagrams in question are then
\begin{itemize}
    \item  same number of dots in each row

    \item  all dots are matched

    \item  no uplines

    \item  the \emph{arcs} in each row are standard:
\begin{pspicture}(0,0)(1,.5)
    \psdots(0,0)(1,0)
\psbezier[linewidth=1pt](0,0)(.3,.3)(.7,.3)(1,0)
\end{pspicture}\quad
\begin{pspicture}(0,0)(2,.5)
    \psdots(0,0)(2,0)
\psbezier[linewidth=.8pt](0,0)(.6,.5)(1.4,.5)(2,0)
\end{pspicture}\quad
\begin{pspicture}(0,0)(1,.5)
    \psdots(0,0)(1,0)
\psbezier[linewidth=1pt](0,0)(.3,.3)(.7,.3)(1,0)
\end{pspicture}\ldots (no crossings or 
    nestings)
    
\end{itemize}
The parameters $n,k,j$ appear respectively as number of dots in each 
row, number of arcs in each row, and the position in the bottom row of 
its last dot incident with an arc. As noted above, the labels are not necessary 
and serve only for identification.

Now we will give a bijection from these dot diagrams to the Cartesian 
product of $\s(j,2k)$, the Stirling matchings defined above, and 
$\p(2k+1,n-j)$, the perfect-power matchings defined above. Since, by Props. \ref{p1} 
and \ref{p2}, $\s(j,2k)$ and 
$\p(2k+1,n-j)$ are counted by $\StirlingPartition{j}{2k}$ and 
$(2k+1)^{n-j}$ respectively, the Theorem will follow.

To get the Stirling matching, take the first $j-1$ dots in each row 
and the \emph{lines} connecting them.

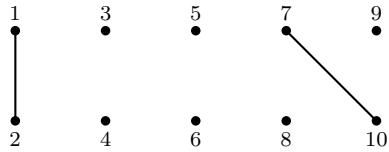
\begin{center} 
    
\begin{pspicture}(-2,-0.8)(2,1.5)
\psset{xunit=1.2cm,yunit=.8cm}

\psline(-2,1)(-2,-0.5)
\psline(1,1)(2,-0.5)

\psdots(-2,1)(-1,1)(0,1)(1,1)(2,1)
\psdots(-2,-0.5)(-1,-0.5)(0,-0.5)(1,-0.5)(2,-0.5)
 
\rput(-2,1.3){\textrm{{\scriptsize $1$}}}
\rput(-1,1.3){\textrm{{\scriptsize $3$}}}
\rput(0,1.3){\textrm{{\scriptsize $5$}}}
\rput(1,1.3){\textrm{{\scriptsize $7$}}}
\rput(2,1.3){\textrm{{\scriptsize $9$}}}

\rput(-2,-0.8){\textrm{{\scriptsize $2$}}}
\rput(-1,-0.8){\textrm{{\scriptsize $4$}}}
\rput(0,-0.8){\textrm{{\scriptsize $6$}}}
\rput(1,-0.8){\textrm{{\scriptsize $8$}}}
\rput(2,-0.8){\textrm{{\scriptsize $10$}}}

\end{pspicture}
\end{center} 
Notice that there may be vertical lines but now eliminate them
by making the technical adjustment of shifting the bottom row one 
unit to the right and adding a dot to each row:

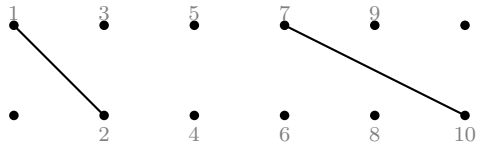
\begin{center} 
    
\begin{pspicture}(-2,-0.5)(3,1.2)
\psset{xunit=1.2cm,yunit=.8cm}

\psline(-2,1)(-1,-0.5)
\psline(1,1)(3,-0.5)

\psdots(-2,1)(-1,1)(0,1)(1,1)(2,1)(3,1)
\psdots(-2,-0.5)(-1,-0.5)(0,-0.5)(1,-0.5)(2,-0.5)(3,-0.5)
\gray{ 
\rput(-2,1.2){\textrm{{\scriptsize $1$}}}
\rput(-1,1.2){\textrm{{\scriptsize $3$}}}
\rput(0,1.2){\textrm{{\scriptsize $5$}}}
\rput(1,1.2){\textrm{{\scriptsize $7$}}}
\rput(2,1.2){\textrm{{\scriptsize $9$}}}
\rput(-1,-0.8){\textrm{{\scriptsize $2$}}}
\rput(0,-0.8){\textrm{{\scriptsize $4$}}}
\rput(1,-0.8){\textrm{{\scriptsize $6$}}}
\rput(2,-0.8){\textrm{{\scriptsize $8$}}}
\rput(3,-0.8){\textrm{{\scriptsize $10$}}} }

\end{pspicture}
\end{center} 
This is the $(j,2k)$ Stirling matching.

To get the perfect-power matching, delete the first $j$ dots in the 
bottom row and the dots in the top row connected to them, and then 
delete all arcs in the top row while leaving their endpoints intact:
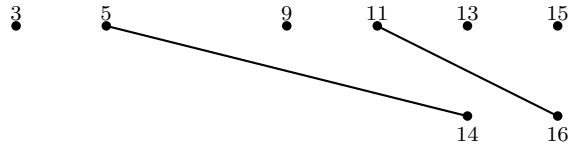
\begin{center} 
    
\begin{pspicture}(-4,-0.7)(4,1.3)
\psset{xunit=1.2cm,yunit=.8cm}

\psline(-2,1)(2,-0.5)

\psline(1,1)(3,-0.5)
\psdots(-3,1)(-2,1)(0,1)(1,1)(2,1)(3,1)
\psdots(2,-0.5)(3,-0.5)

\rput(-3,1.2){\textrm{{\scriptsize $3$}}}
\rput(-2,1.2){\textrm{{\scriptsize $5$}}}

\rput(0,1.2){\textrm{{\scriptsize $9$}}}
\rput(1,1.2){\textrm{{\scriptsize $11$}}}
\rput(2,1.2){\textrm{{\scriptsize $13$}}}
\rput(3,1.2){\textrm{{\scriptsize $15$}}}

\rput(2,-0.8){\textrm{{\scriptsize $14$}}}
\rput(3,-0.8){\textrm{{\scriptsize $16$}}}

\end{pspicture}
\end{center} 

Make the same technical adjustment of shifting the bottom row and 
``prettify'' the diagram:
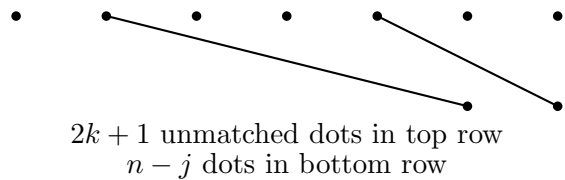
\begin{center} 
    
\begin{pspicture}(-4,-1.2)(4,1.3)
\psset{xunit=1.2cm,yunit=.8cm}
\psline(-2,1)(2,-0.5)

\psline(1,1)(3,-0.5)
\psdots(-3,1)(-2,1)(-1,1)(0,1)(1,1)(2,1)(3,1)
\psdots(2,-0.5)(3,-0.5)

\rput(0,-1){\textrm{{\small $2k+1$ unmatched dots in top row}}}
\rput(0,-1.5){\textrm{{\small $n-j$ dots in bottom row}}}

\end{pspicture}
\end{center} 
This is the $(2k+1,n-j)$ perfect-power matching.

It is easy to check that these maps define a bijection from no-upline 
PM dot diagrams to $\s(j,2k)\times \p(2k+1,n-j)$ and, 
as noted, the Theorem follows. \qed 
 
\section{An Interpretation of Klazar's Recurrence}
In the preceding section we showed that the number $a(n)$ of no-upline 
PM dot diagrams of size $n$ has the same \gf as $w_{12}(n)$. In this 
section we show directly that $a(n)$ satisfies recurrence (\ref{eq2}).
Write (\ref{eq2}) in the equivalent form
\begin{equation}
    a(n)=\underbrace{a(n-1)\vphantom{\binom{n-1}{k+2} }}_{\textrm{\small (1)}} + 
\underbrace{(n-1)a(n-1)\vphantom{\binom{n-1}{k+2} 
}}_{\textrm{\small (2)}} + 
\underbrace{\sum_{k=0}^{n-3}\binom{n-1}{k+2}a(n-2-k)}_{\textrm{\small (3)}}
    \label{eq4}
\end{equation}
and split the no-upline PM dot diagrams of size $n$ into 3 classes according to 
the partners (matched entries) of the last dots in each row, $n$ top
and $n$ bot, as follows. Recall that the partner of $n$ top, denoted 
$p(n$ top), is necessarily in the top 
row unless $n$ top is matched to $n$ bot. Class (1) consists of the dot diagrams in which 
$n$ top is matched to $n$ bot, that is, the last dots in each row 
are joined by a vertical line. Class (2) consists of the dot diagrams in 
which $p(n$ bot) is either (i) in the top row or (ii) in the 
bottom row and $p(n$ top) $<$ $p(n$ bot). 
Class (3) consists of the dot diagrams in 
which $p(n$ top)  $>$ $p(n$ bot) and $p(n$ bot)  is in the 
bottom row. (See the illustrations below.) Now let us count the dot 
diagrams in each class. 

For a dot diagram in Class (1), delete the last dot in each row. This 
gives a dot diagram of size $n-1$---$a(n-1)$ possibilities---and the 
original dot diagram can of course be uniquely recovered from it.

For a dot diagram in Class (2), highlight the partner of $n$ top with a 
heavy dot, then delete the last dot in each row (and the arcs/lines 
therefrom) and join up their partners. The result is a dot diagram of size 
$n-1$ with one highlighted dot in the top row---$(n-1)a(n-1)$ 
possibilities---which uniquely determines the original.
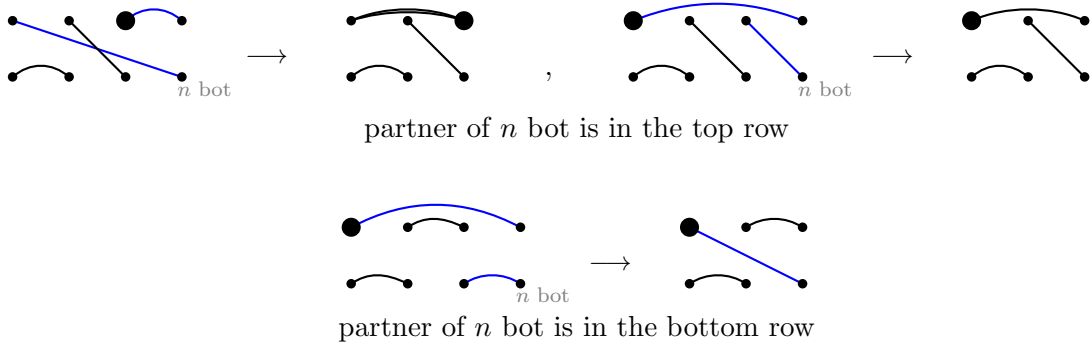
\begin{center} 
    
\begin{pspicture}(-8.4,-3.9)(10.6,.3)
    
\psset{xunit=.75cm,yunit=.5cm}

\psline[linecolor=blue](-10,1)(-7,-0.5)
\psline[linecolor=blue](3,1)(4,-0.5)

\psbezier[linecolor=blue,linewidth=.8pt](1,1)(2,1.6)(3,1.6)(4,1) 
\psbezier[linecolor=blue,linewidth=.8pt](-8,1)(-7.7,1.4)(-7.3,1.4)(-7,1)  
 
\psbezier[linewidth=.8pt](-4,1)(-3.4,1.4)(-2.6,1.4)(-2,1)
\psbezier[linewidth=.8pt](7,1)(7.6,1.4)(8.4,1.4)(9,1)
 
\psbezier[linewidth=.8pt](1,-0.5)(1.3,-.1)(1.7,-.1)(2,-0.5)
\psbezier[linewidth=.8pt](7,-0.5)(7.3,-.1)(7.7,-.1)(8,-0.5)
\psbezier[linewidth=.8pt](-10,-0.5)(-9.7,-.1)(-9.3,-.1)(-9,-0.5)
\psbezier[linewidth=.8pt](-4,-0.5)(-3.7,-.1)(-3.3,-.1)(-3,-0.5)

\psbezier[linewidth=.8pt](-4,1)(-3.4,1.3)(-2.6,1.3)(-2,1)

\psdots(-10,1)(-9,1)(-8,1)(-7,1)(-4,1)(-3,1)(-2,1)(1,1)(2,1)(3,1)(4,1)(7,1)(8,1)(9,1)
\psdots(-10,-0.5)(-9,-0.5)(-8,-0.5)(-7,-0.5)(-4,-0.5)(-3,-0.5)(-2,-0.5)(1,-0.5)(2,-0.5)(3,-0.5)(4,-0.5)(7,-0.5)(8,-0.5)(9,-0.5)

\psline(-9,1)(-8,-0.5)
\psline(-3,1)(-2,-0.5)
\psline(2,1)(3,-0.5)
\psline(8,1)(9,-0.5)
 
\rput(-0.5,-.5){,}

\rput(0,-1.9){\textrm{{\small  partner of $n$ bot is in the top row}}}

\rput(-5.5,0){\textrm{{\footnotesize  $\longrightarrow$}}}
\rput(5.5,0){\textrm{{ \footnotesize  $\longrightarrow$}}}
\rput(0.5,-5.5){\textrm{{ \footnotesize  $\longrightarrow$}}}

\psdots[dotscale=2](-8,1)(-2,1)
\psdots[dotscale=2](1,1)(7,1)

\psline[linecolor=blue](2,-4.5)(4,-6)
\psbezier[linecolor=blue,linewidth=.8pt](-4,-4.5)(-3,-3.7)(-2,-3.7)(-1,-4.5)  
\psbezier[linecolor=blue,linewidth=.8pt](-2,-6)(-1.7,-5.7)(-1.3,-5.7)(-1,-6)

\psbezier[linewidth=.8pt](-3,-4.5)(-2.7,-4.2)(-2.4,-4.2)(-2,-4.5)
\psbezier[linewidth=.8pt](-4,-6)(-3.7,-5.7)(-3.4,-5.7)(-3,-6)

\psbezier[linewidth=.8pt](3,-4.5)(3.3,-4.2)(3.7,-4.2)(4,-4.5)
\psbezier[linewidth=.8pt](2,-6)(2.3,-5.7)(2.7,-5.7)(3,-6)
 
\psdots(-4,-4.5)(-3,-4.5)(-2,-4.5)(-1,-4.5)(2,-4.5)(3,-4.5)(4,-4.5)
\psdots(-4,-6)(-3,-6)(-2,-6)(-1,-6)(2,-6)(3,-6)(4,-6)
\psdots[dotscale=2](-4,-4.5)(2,-4.5)

\rput(0,-7.2){\textrm{{\small  partner of $n$ bot is in the bottom row}}}

{\gray
\rput(-6.6,-.8){\textrm{{\scriptsize $n$ bot}}}
\rput(4.4,-.8){\textrm{{\scriptsize $n$ bot}}}
\rput(-0.6,-6.3){\textrm{{\scriptsize $n$ bot}}}}

\end{pspicture}
\end{center}

For a dot diagram in Class (3), first record the locations of the partners of 
the last dots and of all $k \ge 0$ \emph{vertical} lines lying \emph{between} these partners. 
This gives a $(k+2)$-element subset $X$ of $[n-1]$. Then delete the last 
dots, their partners, all these vertical lines, and ``prettify'' the 
diagram. The result is a dot diagram of size $n-k-2$ which, together 
with the set $X$, determines the original. 
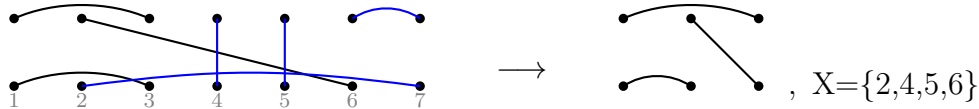
\begin{center} 
    
\begin{pspicture}(-7,-0.7)(7,1)
\psset{xunit=.9cm,yunit=.6cm}

\psline(-6,1)(-2,-0.5)
\psline(3,1)(4,-0.5)
  
\psbezier[linewidth=.8pt](-7,1)(-6.4,1.4)(-5.6,1.4)(-5,1)
\psbezier[linewidth=.8pt](-7,-0.5)(-6.4,-.1)(-5.6,-.1)(-5,-0.5) 

\psbezier[linewidth=.8pt](2,1)(2.6,1.4)(3.4,1.4)(4,1) 
\psbezier[linewidth=.8pt](2,-0.5)(2.3,-.2)(2.7,-.2)(3,-0.5)   

\psdots(-7,1)(-6,1)(-5,1)(-4,1)(-3,1)(-2,1)(-1,1)(2,1)(3,1)(4,1)
\psdots(-7,-0.5)(-6,-0.5)(-5,-0.5)(-4,-0.5)(-3,-0.5)(-2,-0.5)(-1,-0.5)(2,-0.5)(3,-0.5)(4,-0.5)

\psline[linecolor=blue](-4,1)(-4,-0.5)
\psline[linecolor=blue](-3,1)(-3,-0.5)
\psbezier[linecolor=blue,linewidth=.8pt](-6,-0.5)(-4,-0.1)(-3,-0.1)(-1,-0.5) 
\psbezier[linecolor=blue,linewidth=.8pt](-2,1)(-1.7,1.3)(-1.3,1.3)(-1,1)  

\rput(0.5,-0.2){$\longrightarrow$}

\rput(6.1,-0.5){X=\{2,4,5,6\} }
\rput(4.5,-0.7){,}

\gray{
\rput(-7,-0.8){\textrm{{\scriptsize $1$}}}
\rput(-6,-0.8){\textrm{{\scriptsize $2$}}}
\rput(-5,-0.8){\textrm{{\scriptsize $3$}}}
\rput(-4,-0.8){\textrm{{\scriptsize $4$}}}
\rput(-3,-0.8){\textrm{{\scriptsize $5$}}}
\rput(-2,-0.8){\textrm{{\scriptsize $6$}}}
\rput(-1,-0.8){\textrm{{\scriptsize $7$}}}
 }

\end{pspicture}
\end{center} 
Notice that in this case the vertical lines have to be deleted, for 
otherwise the resulting dot diagram would contain lines of positive 
slope, which are forbidden. 

Thus the three classes are counted by the three terms in 
(\ref{eq4}). \qed

\section{Codes for Trees and Matchings}
An \o of size $n$ can be built up from the root 0 by successively 
adding vertices $1,2,\ldots,n$. Vertex 1 is necessarily a (rightmost) 
child of the root and for $2\le i \le n$, vertex $i$ is either the 
rightmost child of a vertex $v \in [0,i-1]$---coded as $(R,v)$---or the 
(immediate) left neighbor of a vertex $v \in [1,i-1]$---coded as 
$(L,v)$. Thus an \o of size $n$ corresponds naturally to a 
\emph{build-tree code} $\big( (X_{k},i_{k}) \big)_{1\le k \le 
n}$ where $(X_{1},i_{1})=(R,0)$  and for $2\le k \le n$, $X_{k}=R$ 
and $i_{k}\in[0,k-1]$ or $X_{k}=L$ 
and $i_{k}\in[1,k-1]$. For example, with $X_{v}$ short for $(X,v)$,
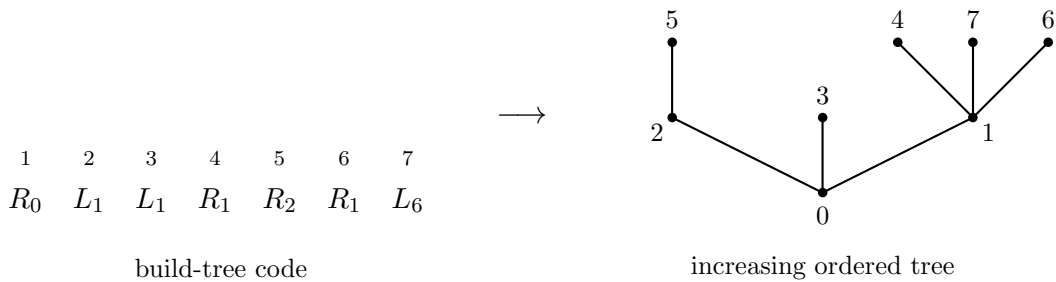
\begin{center} 
\setlength{\arraycolsep}{2mm}
\begin{pspicture}(-2,-2)(15.5,2.5)
\psline(8,2)(8,1)(10,0)(10,1)
\psline(10,0)(12,1)(11,2)
\psline(12,2)(12,1)(13,2)

\psdots(8,1)(8,2)(10,0)(10,1)(11,2)(12,1)(12,2)(13,2)

\rput(10,-0.3){\textrm{{\footnotesize $0$}}}
\rput(12.2,0.8){\textrm{{\footnotesize $1$}}}
\rput(7.8,0.8){\textrm{{\footnotesize $2$}}}
\rput(10,1.3){\textrm{{\footnotesize $3$}}}
\rput(11,2.3){\textrm{{\footnotesize $4$}}}
\rput(8,2.3){\textrm{{\footnotesize $5$}}}
\rput(13,2.3){\textrm{{\footnotesize $6$}}}
\rput(12,2.3){\textrm{{\footnotesize $7$}}}

\rput(2,0.2){\textrm{{\small $
\begin{array}{ccccccc}
    \textrm{{\scriptsize 1}} &  \textrm{{\scriptsize 2}} &  \textrm{{\scriptsize 
    3}} &  \textrm{{\scriptsize 4}} &  \textrm{{\scriptsize 5}} &  \textrm{{\scriptsize 6}} &  \textrm{{\scriptsize 7}} \\
    R_{0} & L_{1} &  L_{1} & R_{1} & R_{2} & R_{1} &  L_{6}
\end{array}
$
}}}

\rput(6,1){$\longrightarrow$}

\rput(10,-1){\textrm{{\footnotesize \o}}}
\rput(2,-1){\textrm{{\footnotesize build-tree code}}}

\rput(5.8,-2){\textrm{{\small natural correspondence between build-tree 
codes and \os}}}

\end{pspicture}
\end{center}

To reverse this correspondence, if $n$ is a rightmost child, record 
$(X_{n},i_{n})=(R,i)$ with $i$ the parent of $n$, otherwise 
$(X_{n},i_{n})=(L,i)$ with $i$ the right neighbor of $n$. Delete $n$ 
and proceed similarly, starting with $n-1$, to obtain $(X_{k},i_{k})$ 
for $k=n,n-1,\ldots,1$ in turn.

Similarly, a PM dot diagram can be built up from the empty diagram 
by successively adding a rightmost dot to each row, connecting the new 
top dot to one of the dots $i$ in the bottom row---coded as 
$(B,i)$---or to one of the old dots $i$ in the top row---coded as 
$(T,i)$. Thus a diagram of size $n$ corresponds naturally to a 
\emph{build-matching code} $\big( (Y_{k},i_{k}) \big)_{1\le k \le 
n}$ where $(Y_{1},i_{1})=(B,1)$  and for $2\le k \le n$, $Y_{k}=B$ 
and $i_{k}\in[1,k]$ or $Y_{k}=T$ and $i_{k}\in[1,k-1]$. 

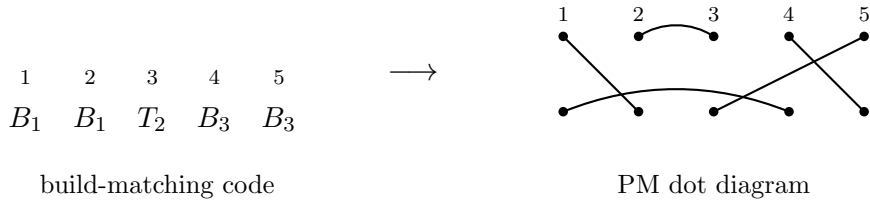
\begin{center} 
\setlength{\arraycolsep}{2mm}
\begin{pspicture}(-8,-2.1)(8,1)
\psline(2,1)(3,0)
\psline(4,0)(6,1)
\psline(5,1)(6,0)

\psbezier[linewidth=.8pt](3,1)(3.3,1.2)(3.7,1.2)(4,1)
\psbezier[linewidth=.8pt](2,0)(3,0.4)(4,0.4)(5,0)

\psdots(2,0)(3,0)(4,0)(5,0)(6,0)
\psdots(2,1)(3,1)(4,1)(5,1)(6,1)

\rput(2,1.3){\textrm{{\scriptsize $1$}}}
\rput(3,1.3){\textrm{{\scriptsize $2$}}}
\rput(4,1.3){\textrm{{\scriptsize $3$}}}
\rput(5,1.3){\textrm{{\scriptsize $4$}}}
\rput(6,1.3){\textrm{{\scriptsize $5$}}}

\rput(-3.4,0.2){\textrm{{\small $
\begin{array}{ccccc}
    \textrm{{\scriptsize 1}} &  \textrm{{\scriptsize 2}} &  \textrm{{\scriptsize 
    3}} &  \textrm{{\scriptsize 4}} &  \textrm{{\scriptsize 5}}  \\
    B_{1} & B_{1} &  T_{2} & B_{3} & B_{3}
\end{array}
$
}}}

\rput(0,0.5){$\longrightarrow$}

\rput(4,-1){\textrm{{\footnotesize PM dot diagram}}}
\rput(-3.4,-1){\textrm{{\footnotesize build-matching code}}}

\rput(0,-2){\textrm{{\small natural correspondence between build-matching 
codes and PM dot diagrams}}}

\end{pspicture}
\end{center}

Since there are $2k-1$ possibilities for $(X_{k},i_{k})$ 
in a  build-tree code and for $(Y_{k},i_{k})$ in a build-matching 
code, both are counted by $(2n-1)!!$, and we see once again that there 
are $(2n-1)!!$ \os of size $n$. 

We need the correspondence between the two codes obtained by 
identifying $R \leftrightarrow B$ and $L \leftrightarrow  T$ 
and, when $(X_{k},i_{k})$ in a build-tree code has $X_{k}=R$ and
$i_{k}=0$, replacing it by $(Y_{k},i_{k})=(B,k)$.
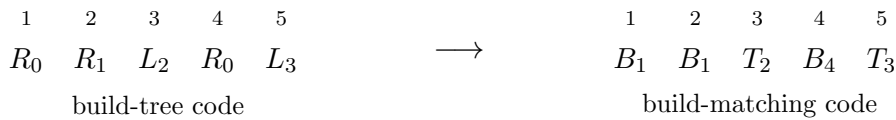
\begin{center} 
\setlength{\arraycolsep}{2mm}
\begin{pspicture}(-8,-1.7)(8,1)

\rput(-4,0.2){\textrm{{\small $
\begin{array}{ccccc}
    \textrm{{\scriptsize 1}} &  \textrm{{\scriptsize 2}} &  \textrm{{\scriptsize 
    3}} &  \textrm{{\scriptsize 4}} &  \textrm{{\scriptsize 5}}  \\
    R_{0} & R_{1} &  L_{2} & R_{0} & L_{3}
\end{array}
$
}}}

\rput(4,0.2){\textrm{{\small $
\begin{array}{ccccc}
    \textrm{{\scriptsize 1}} &  \textrm{{\scriptsize 2}} &  \textrm{{\scriptsize 
    3}} &  \textrm{{\scriptsize 4}} &  \textrm{{\scriptsize 5}}  \\
    B_{1} & B_{1} &  T_{2} & B_{4} & T_{3}
\end{array}
$
}}}

\rput(0,0){$\longrightarrow$}

\rput(4,-0.7){\textrm{{\footnotesize build-matching code}}}
\rput(-4,-0.7){\textrm{{\footnotesize build-tree code}}}

\rput(0,-1.5){\textrm{{\small correspondence between build-tree
codes and build-matching codes}}}

\end{pspicture}
\end{center}

\section{A Combinatorial Proof of Klazar's Recurrence}
\subsection{Preliminaries}
To complete the combinatorial proof of recurrence (\ref{eq2}) in 
context we need
a bijection from Klazar trees to perfect 
matchings with no even-to-odd matches, equivalently to PM dot diagrams 
with no uplines. In fact, we will give a 
bijection $\Phi$ from \os to PM dot diagrams that sends Klazar violators 
to uplines. More precisely, define the \emph{partner} of 
a Klazar violator to be its rightmost child or closest left sibling, 
whichever is larger if both are present (for a Klazar violator at 
least one must be present). For example, the  tree below has 4 Klazar 
violators with partners as shown.

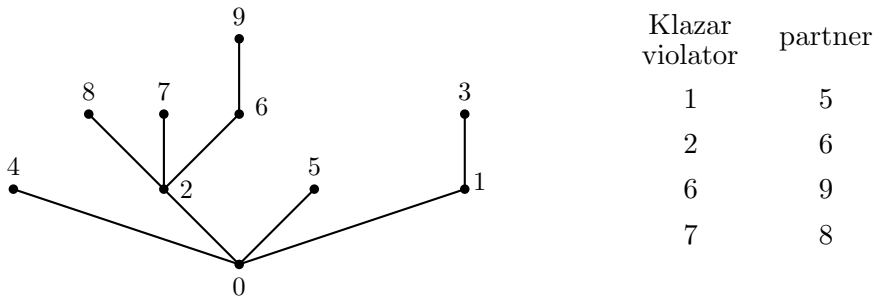
\begin{center} 

\begin{pspicture}(-3,0)(9,3.5)
\psline(-3,1)(0,0)(3,1)(3,2)
\psline(-2,2)(-1,1)(0,0)(1,1)
\psline(-1,2)(-1,1)(0,2)(0,3)

\psdots(-3,1)(-2,2)(-1,1)(-1,2)(0,2)(0,3)(1,1)(3,1)(0,0)(3,2)

\rput(0,-0.3){\textrm{{\footnotesize $0$}}}
\rput(3.2,1.1){\textrm{{\footnotesize $1$}}}
\rput(-0.7,1,2.3){\textrm{{\footnotesize $2$}}}
\rput(3,2.3){\textrm{{\footnotesize $3$}}}
\rput(-3,1.3){\textrm{{\footnotesize $4$}}}
\rput(1,1.3){\textrm{{\footnotesize $5$}}}
\rput(0.3,2.1){\textrm{{\footnotesize $6$}}}
\rput(-1,2.3){\textrm{{\footnotesize $7$}}}
\rput(-2,2.3){\textrm{{\footnotesize $8$}}}
\rput(0,3.3){\textrm{{\footnotesize $9$}}}

\rput(6,3.2){\textrm{{\small Klazar}}}
\rput(6,2.8){\textrm{{\small violator}}}

\rput(6,2.2){\textrm{{\small 1}}}
\rput(6,1.6){\textrm{{\small 2}}}
\rput(6,1.0){\textrm{{\small 6}}}
\rput(6,0.4){\textrm{{\small 7}}}

\rput(7.8,3.0){\textrm{{\small partner}}}

\rput(7.8,2.2){\textrm{{\small 5}}}
\rput(7.8,1.6){\textrm{{\small 6}}}
\rput(7.8,1.0){\textrm{{\small 9}}}
\rput(7.8,0.4){\textrm{{\small 8}}}

\end{pspicture}
\end{center}

Then the bijection $\Phi$ sends each Klazar violator $i$ together with its partner $j$ to 
an upline from the $i$th dot in the bottom row to the $j$th dot in the 
top row. 

\begin{prop}
    In an \o, the map Klazar violator $\mapsto$ partner is one-to-one.
\end{prop}
Proof.  A vertex in an \o is the partner of at most one Klazar violator 
because a partner is either a rightmost child or a left neighbor and a 
vertex cannot simultaneously be both a  rightmost child and a left 
neighbor. \qed

A \emph{child} vertex in an \o is a vertex that is a child of some 
other vertex, that is, a non-root vertex.
Partition the child vertices in an \o into two classes: those that 
are the partner of some Klazar violator (``partners'') and those that 
are not (``non-partners''). 
\begin{prop}
   In an \o, there is a bijection $H$ from Klazar compliant child vertices to non-partners.
   \label{H}
\end{prop}
Proof. Suppose $v$ is a  Klazar compliant child vertex. If $v$ is a 
non-partner, set $H(v)=v$. Otherwise, $v_{1}:=v$ is the partner of a (unique) 
Klazar violator $v_{2}$ and clearly $v_{2}<v_{1}$. If $v_{2}$ is a 
non-partner, set $H(v)=v_{2}$. Otherwise proceed similarly to obtain 
vertices
$v_{2}>v_{3}> \ldots > v_{k}$ stopping at the first $v_{k}$ that is a 
non-partner and set 
$H(v)=v_{k}$. For example, the tree above has 5 Klazar compliant 
child vertices and 5 non-partners as in the table:
\begin{center}
    \begin{tabular}{cccccc}    
        Klazar compliant child  $v$ & 3 & 4 & 5 & 8 & 9  \\     
        non-partner  $H(v)$ & 3 & 4 & 1 & 7 & 2  \\
    \end{tabular}
\end{center}
\qed

\begin{prop}
  There is an involution $F$ on \os of size $n$ that is the identity 
on trees in which $n$ does not have a right neighbor and 
otherwise flips the Klazar violator/complier status of the right neighbor 
of $n$ while not disturbing any other Klazar violators or their 
partners.
\label{F}
\end{prop}
Proof. In case the right neighbor $j$ of $n$ is a Klazar violator 
transfer the associate of 
$j$ and its big cohort so that they become the leftmost segment of the children of $j$; in 
case the right neighbor $j$ of $n$ is a Klazar complier transfer the 
smallest child of $j$ and its cohort to the 
immediate left of the big cohort of $j$, as illustrated.

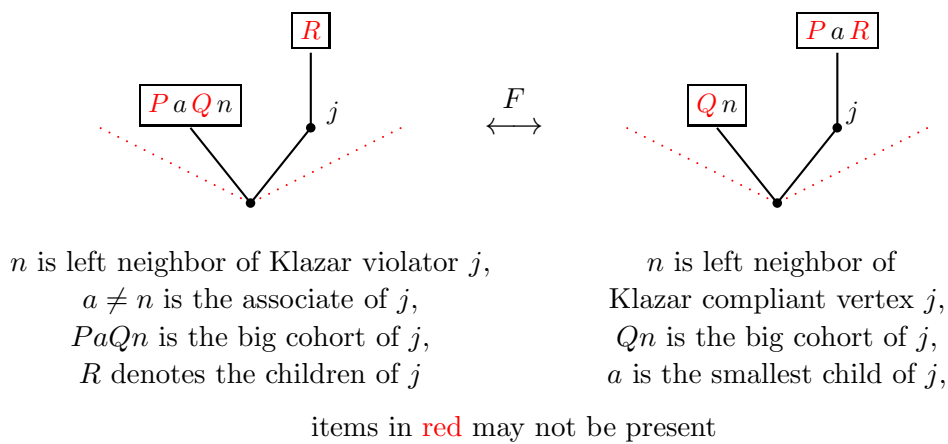
\begin{center} 

\begin{pspicture}(-6.5,-4)(6.5,2.5)
\psline[linestyle=dotted,linecolor=red](-3.5,0)(-5.5,1)
\psline[linestyle=dotted,linecolor=red](-3.5,0)(-1.5,1)
\psline[linestyle=dotted,linecolor=red](3.5,0)(1.5,1)
\psline[linestyle=dotted,linecolor=red](3.5,0)(5.5,1)
\psline(-3.5,0)(-4.3,1)
\psline(-3.5,0)(-2.7,1)(-2.7,2)
\psline(3.5,0)(2.7,1)
\psline(3.5,0)(4.3,1)(4.3,2)

\psdots(-3.5,0)(-2.7,1)(3.5,0)(4.3,1)

\rput(-2.4,1.2){\textrm{{\footnotesize $j$}}}
\rput(-2.7,2.3){\textrm{{\footnotesize \framebox{$\red{R}$}}}}
\rput(-4.3,1.3){\textrm{{\footnotesize \framebox{$\red{P}\,a\,\red{Q}\,n$}}}}
\rput(2.7,1.3){\textrm{{\footnotesize \framebox{$\red{Q}\,n$}}}}
\rput(4.3,2.3){\textrm{{\footnotesize \framebox{$\red{P}\,a\,\red{R}$}}}}
\rput(4.5,1.2){\textrm{{\footnotesize $j$}}}

\rput(0,1){$\longleftrightarrow$}
\rput(0,1.4){\textrm{{\small $F$}}}

\rput(-3.5,-0.8){\textrm{{\small $n$ is left neighbor of Klazar 
violator $j$,}}}
\rput(-3.5,-1.3){\textrm{{\small $a \ne n$ is the associate of $j$,}}}
\rput(-3.5,-1.8){\textrm{{\small $PaQn$ is the big cohort of $j$,}}}
\rput(-3.5,-2.3){\textrm{{\small $R$ denotes the children of $j$}}}

\rput(3.5,-0.8){\textrm{{\small $n$ is left neighbor of }}}
\rput(3.5,-1.3){\textrm{{\small Klazar compliant vertex $j$,}}}
\rput(3.5,-1.8){\textrm{{\small $Qn$ is the big cohort of $j$,}}}
\rput(3.5,-2.3){\textrm{{\small $a$ is the smallest child of $j$,}}}

\rput(0,-3){\textrm{{\small items in \red{red} may not be present}}}
\rput(0,-4){\textrm{ The involution $F$}}

\end{pspicture}
\end{center} 
\qed

To prune an \o $T$ of size $n$ means to delete $n$ (necessarily a 
leaf) and its incident edge. We use $P(T)$ to denote the pruned tree. 
Also, for a Klazar violator $v$, $p(v)$ denotes its partner.

A PM dot diagram of size $n-1$ can be enlarged to one of size $n$ in $2n-1$ ways:
add a dot at the end of each row and then either join these two dots 
together 
leaving the rest of the original \pm intact, or join the new dot in the 
top row to any one of the $2n-2$ original dots $i$, delete the 
edge from $i$ to its original partner $j$ and replace it with an edge from $j$ to the 
new dot in the bottom row.
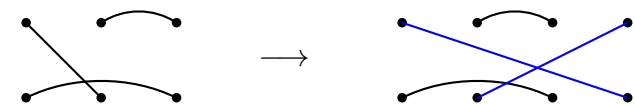
\begin{center} 
    \label{enlarge}
    \begin{pspicture}(-4,-1)(4,1.3)
    
\psbezier[linewidth=.8pt](-3,1)(-2.7,1.2)(-2.3,1.2)(-2,1)
\psbezier[linewidth=.8pt](-4,0)(-3.4,0.3)(-2.6,0.3)(-2,0)
\psbezier[linewidth=.8pt](2,1)(2.3,1.2)(2.7,1.2)(3,1)
\psbezier[linewidth=.8pt](1,0)(1.6,0.3)(2.4,0.3)(3,0)

\psdots(-4,1)(-4,0)(-3,1)(-3,0)(-2,1)(-2,0)
\psdots(1,0)(1,1)(2,0)(2,1)(3,0)(3,1)(4,0)(4,1)

\psline(-4,1)(-3,0)
\psline[linecolor=blue](1,1)(4,0)
\psline[linecolor=blue](4,1)(2,0)

\rput(-.5,.5){$\longrightarrow$ }
 
\rput(0,-0.9){\small{enlarge using ``2 bot'' (second dot in bottom row)} }

\end{pspicture}
\end{center} 

\subsection{Recursive $\Phi$}
Now we can give a recursive definition of $\Phi$. First, $\Phi$ sends 
the unique \o of one edge to the unique PM dot diagram of size 1. For an \o $T$ of 
size $n \ge 2$, define $\Phi(T)$ according to the 6 cases in the 
following table. Recall that $P$ (for prune), $F$ (an involution), $p$ (for 
partner), and the  bijection $H$ have all been defined in the 
preceding subsection and $KV,\ KC$ mean Klazar violator, complier respectively.
\begin{center}
    \begin{tabular}{|c|cccc|}
        \hline
        \raisebox{-.3ex}[0pt]{\#} & \raisebox{-.3ex}[0pt]{Case} & \raisebox{-.3ex}[0pt]{Status of $j$} & 
	\raisebox{-.3ex}[0pt]{Enlarge:} & \raisebox{-.3ex}[0pt]{using:}  \\[2mm]
        \hline
        1 & $n$ rightmost child of root & - & $\Phi(P(T))$ & $n$ 
        bot  \\[1mm]  \hline
        2 & $n$ rightmost child of KV $j$ & $j$ is KV  in $P(T)$ & $\Phi(P(T))$
        &   $j$ bot  \\[1mm]  \hline
        & \raisebox{-.3ex}[0pt]{$n$ is left neighbor of KV $j$} &  & 
         &   \\
	\raisebox{1.5ex}[0pt]{3}  & and $n$ is associate of $j$ & 
	\raisebox{1.5ex}[0pt]{$j$ is KC leaf in $P(T)$} &  
	\raisebox{1.5ex}[0pt]{$\Phi(P(T))$} &  \raisebox{1.5ex}[0pt]{$j$ bot} \\[1mm] \hline
	
        &   & \raisebox{-.3ex}[0pt]{$j$ is KC and $H(j)$ is} &  &    \\
	
	\raisebox{1.5ex}[0pt]{4} & \raisebox{1.5ex}[0pt]{$n$ rightmost child of 
	KC $j$} &  non-partner in $P(T)$ & 
        \raisebox{1.5ex}[0pt]{$\Phi(P(T))$} & \raisebox{1.5ex}[0pt]{$H(j)$ top}  \\[1mm] \hline
        & \raisebox{-.3ex}[0pt]{$n$ is left neighbor of KV $j$} & 
        \raisebox{-.3ex}[0pt]{$j$ is KC non-leaf}  &  &   \\
         \raisebox{1.5ex}[0pt]{5} & and $n$ is not associate of $j$  & in 
        $P(F(T))$ & \raisebox{1.5ex}[0pt]{$\Phi(P(F(T))$}  & 
        \raisebox{1.5ex}[0pt]{$j$ bot} \\[1mm] \hline
	 6 & $n$ is left neighbor of KC $j$ & $j$ is KV in 
        $P(F(T))$ & $\Phi(P(F(T))$ & $p(j)$ top  \\[1mm]\hline
    \end{tabular}
\end{center}

\centerline{The bijection $\Phi$}
\vspace*{2mm}

\begin{theorem}
    $\Phi$ is a size-preserving bijection from \os to perfect matchings 
    that sends $($Klazar violator, partner$\,)$ pairs $(i,j)$ to uplines 
    $i\nearrow j$.
\end{theorem}
\textbf{Proof}\quad To show that $\Phi$ is a bijection it suffices, by induction, to show 
that for a given \o $T$ of size $n$, the specifications in the table above will 
enlarge $\Phi(P(T))$ using a full complement of dots: $i$ top, $1\le i \le 
n-1$ and $i$ bot, $1\le i\le n$. In the first four cases $\Phi(P(T))$ 
gets enlarged from $T$ itself and in the last two cases from $F(T)$ 
(recall $F$ is an involution).

Partition the child (non-root) vertices $j$ of $P(T)$ into three classes: (1) 
Klazar violators, (2) Klazar compliant leaves, (3) Klazar compliant 
non-leaves. Cases 2,\,3,\,5 hit $j$ bot for $j$ in classes (1), 
(2), (3) respectively. Now partition the child vertices of $P(T)$
in another way
into two classes: (1$'$) partners (of some Klazar violator), (2$'$) 
non-partners. Case 6 hits $j$ top for all $j$ in class (1$'$), and, 
using Prop. \ref{H}, case 4 hits $j$ top for all $j$ in class (2$'$). 
Since $n$ bot is hit by case 1, this proves that $\Phi$ is a bijection.

Finally, to verify that Klazar violator/partner pairs $(i,j)$ are 
transformed into uplines from $i$ bot to $j$ top is a matter of 
checking cases. For example, in case 6, by induction and Prop. 
\ref{F}, $\Phi(P(F(T)))$ has an  
upline from each Klazar violator to its partner in $P(T)$ (and hence in 
$T$) but also one from $j$ to 
$p(j)$. The latter upline, however, is destroyed when $p(j)$ top is 
used to enlarge $\Phi(P(F(T)))$. Other cases are left to the reader. \qed

\subsection{Explicit $\Phi$}
Here we give an explicit description of $\Phi$ as a composition of 
three bijections: (1) a tweaked version $\si$ of the natural 
bijection from \os to build-tree codes, (2) the correspondence to 
build-matching codes, and (3) a tweaking $\tau$ of the natural 
bijection from build-matching codes to PM dot diagrams.

Definition of \si: Given an \o, obtain a build-tree code just as in 
the natural correspondence except that the involution $F$ is applied 
to the original tree before recording $(X_{n},i_{n})$ and to each 
succeeding pruned tree before recording $(X_{k},i_{k})$. 
Conversely, $\si^{-1}$ builds up the tree in the natural way except 
that $F$ is applied to each intermediate tree before the next one is 
constructed. For example, the last tree in Fig. 3 corresponds under 
\si to the code
$ R_{1},R_{1},L_{1},L_{2},L_{1},R_{2}$ as illustrated ($F$ is not shown
when it is the identity).
\begin{center}

\begin{pspicture}(-8.8,-0.2)(14.2,1.5)

\psset{xunit=.65cm,yunit=.65cm}

\psdots(-13,0)(-10,0)(-10,1)(-7,0)(-7,1)(-7,2)(-4,1)(-3,0)(-2,1)(-2,2)(1,1)(2,0)(2,1)(3,1)(6,1)(7,1)(7.5,0)(8,1)(9,1)

\psline(-10,0)(-10,1)
\psline(-7,0)(-7,1)(-7,2)
\psline(-4,1)(-3,0)(-2,1)(-2,2)
\psline(1,1)(2,0)(3,1)
\psline(2,0)(2,1)
\psline(7,1)(7.5,0)(8,1)
\psline(6,1)(7.5,0)(9,1)

\rput(-11.5,0.5){\textrm{{\footnotesize  $\xrightarrow[R_{0}]{}$}}}
\rput(-8.5,0.5){\textrm{{\footnotesize  $\xrightarrow[R_{1}]{}$}}}
\rput(-5.5,0.5){\textrm{{\footnotesize  $\xrightarrow[L_{1}]{}$}}}
\rput(-0.5,0.5){\textrm{{\footnotesize  $\xrightarrow[\, F\,]{}$}}}
\rput(4.5,0.5){\textrm{{\footnotesize  $\xrightarrow[L_{2}]{}$}}}

\rput(-13,-0.4){\textrm{{\scriptsize 0}}}

\rput(-10,-0.4){\textrm{{\scriptsize 0}}}
\rput(-10,1.4){\textrm{{\scriptsize 1}}}

\rput(-7,-.4){\textrm{{\scriptsize 0}}}
\rput(-6.8,1){\textrm{{\scriptsize 1}}}
\rput(-7,2.4){\textrm{{\scriptsize 2}}}

\rput(-4,1.4){\textrm{{\scriptsize 3}}}
\rput(-3,-.4){\textrm{{\scriptsize 0}}}
\rput(-1.8,1){\textrm{{\scriptsize 1}}}
\rput(-2,2.4){\textrm{{\scriptsize 2}}}

\rput(1,1.4){\textrm{{\scriptsize 2}}}
\rput(2,1.4){\textrm{{\scriptsize 3}}}
\rput(3,1.4){\textrm{{\scriptsize 1}}}
\rput(2,-.4){\textrm{{\scriptsize 0}}}

\rput(6,1.4){\textrm{{\scriptsize 4}}}
\rput(7,1.4){\textrm{{\scriptsize 2}}}
\rput(8,1.4){\textrm{{\scriptsize 3}}}
\rput(9,1.4){\textrm{{\scriptsize 1}}}
\rput(7.5,-.4){\textrm{{\scriptsize 0}}}

\end{pspicture}
\end{center}
\begin{center}

\begin{pspicture}(-8,-1.4)(12,2)

\psset{xunit=.65cm,yunit=.65cm}

\psdots(-6,1)(-5,1)(-4,1)(-3,1)(-2,1)(-4,0)(1,1)(2,0)(2,1)(3,1)(2.5,2)(3.5,2)(6,1)(7,1)(7,0)(8,1)(7.5,2)(8.5,2)(8.5,3)

\psline(-6,1)(-4,0)(-2,1)
\psline(-5,1)(-4,0)(-3,1)
\psline(-4,1)(-4,0)

\psline(1,1)(2,0)(3,1)(3.5,2)
\psline(2,0)(2,1)
\psline(2.5,2)(3,1)

\psline(6,1)(7,0)(7,1)
\psline(7,0)(8,1)(8.5,2)(8.5,3)
\psline(7.5,2)(8,1)

\rput(-7.5,0.5){\textrm{{\footnotesize  $\xrightarrow[L_{1}]{}$}}}
\rput(-0.5,0.5){\textrm{{\footnotesize  $\xrightarrow[\, F\,]{}$}}}
\rput(4.5,0.5){\textrm{{\footnotesize  $\xrightarrow[R_{2}]{}$}}}

\rput(-6,1.4){\textrm{{\scriptsize 4}}}
\rput(-5,1.4){\textrm{{\scriptsize 2}}}
\rput(-4,1.4){\textrm{{\scriptsize 3}}}
\rput(-3,1.4){\textrm{{\scriptsize 5}}}
\rput(-2,1.4){\textrm{{\scriptsize 1}}}
\rput(-4,-.4){\textrm{{\scriptsize 0}}}

\rput(1,1.4){\textrm{{\scriptsize 3}}}
\rput(2,1.4){\textrm{{\scriptsize 5}}}
\rput(3.3,1){\textrm{{\scriptsize 1}}}
\rput(2.5,2.4){\textrm{{\scriptsize 4}}}
\rput(3.5,2.4){\textrm{{\scriptsize 2}}}
\rput(2,-.4){\textrm{{\scriptsize 0}}}

\rput(6,1.4){\textrm{{\scriptsize 3}}}
\rput(7,1.4){\textrm{{\scriptsize 5}}}
\rput(8.3,1){\textrm{{\scriptsize 1}}}
\rput(7.5,2.4){\textrm{{\scriptsize 4}}}
\rput(8.7,2.3){\textrm{{\scriptsize 2}}}
\rput(8.5,3.4){\textrm{{\scriptsize 6}}}
\rput(7,-.4){\textrm{{\scriptsize 0}}}

\rput(0,-2){\textrm{{\small Figure 3}}}

\end{pspicture}
\end{center}
\begin{prop}
    \label{sig}
    Under $\si^{-1}$, the pair $(i,j)$ is a Klazar violator/partner 
    pair iff $(L,i)$ occurs an odd number of times in the build-tree 
    code and $j$ is the position in the code of the last occurrence 
    of an $(X,i)\ (X=L$ or $R)$.
\end{prop}
\emph{Proof}\quad  First, adding a rightmost child from an $(R,-)$ entry in 
the code never introduces a Klazar violator in the \o. On 
the other hand, the 
first occurrence of an $(L,i)$, say as the $j$th entry in the code, makes 
$i$ a Klazar violator with partner $j$. Subsequent occurrences of $(L,i)$ flip the 
violator/compliant status of $i$ while occurrences of either $(L,i)$ or 
$(R,i)$ update the partner of $i$ to the current position in the code 
whenever $i$ is currently a violator. The result follows. \qed

Definition of $\tau$: Given a build-matching code of length $n$, for 
$k=1,2,\ldots,n$ 
enlarge (\ref{enlarge}) the current PM dot diagram (initially empty) using 
the vertex specified as follows. First, suppose $Y_{k}=B$. If 
$i_{k}=k$, use $k$ bot; if $i:=i_{k}<k$ and there is an upline from 
$i$ in the current PM dot diagram, use $i$ bot, otherwise use $j_{1}$ top where 
$j_{1}\nearrow j_{2}\nearrow \ldots \nearrow i$ is the maximal run of 
uplines terminating at $i$ with $j_{1}=i$ if there is no upline terminating 
at $i$. Second, suppose $Y_{k}=T$. If there is an upline $i_{k} \nearrow 
j$  in the current PM dot diagram, use $j$ top, otherwise use $i_{k}$ bot. For 
example, the code $B_{1},T_{1},B_{2},T_{1}$ successively yields 
($\epsilon$ denotes the empty PM dot diagram)
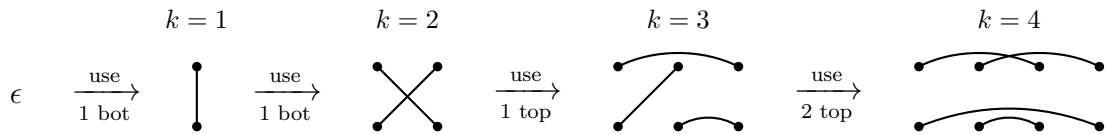
\begin{center} 
\begin{pspicture}(-8,0)(10,1.3)
\psset{unit=.8cm}
\psline(-6,0)(-6,1)
\psline(-3,0)(-2,1)
\psline(-2,0)(-3,1)
\psline(2,1)(1,0)

\psbezier[linewidth=.8pt](1,1)(1.6,1.3)(2.4,1.3)(3,1)
\psbezier[linewidth=.8pt](6,1)(6.6,1.3)(7.4,1.3)(8,1)
\psbezier[linewidth=.8pt](7,1)(7.6,1.3)(8.4,1.3)(9,1)

\psbezier[linewidth=.8pt](6,0)(7,0.4)(8,0.4)(9,0)

\psbezier[linewidth=.8pt](2,0)(2.3,0.2)(2.7,0.2)(3,0)
\psbezier[linewidth=.8pt](7,0)(7.3,0.2)(7.7,0.2)(8,0)

\psdots(-6,0)(-3,0)(-2,0)(1,0)(2,0)(3,0)(6,0)(7,0)(8,0)(9,0)
\psdots(-6,1)(-3,1)(-2,1)(1,1)(2,1)(3,1)(6,1)(7,1)(8,1)(9,1)

\rput(-6,1.8){\textrm{{\footnotesize $k=1$}}}
\rput(-2.5,1.8){\textrm{{\footnotesize $k=2$}}}
\rput(2,1.8){\textrm{{\footnotesize $k=3$}}}
\rput(7.5,1.8){\textrm{{\footnotesize $k=4$}}}

\rput(-9,0.5){$\epsilon$}

\rput(-7.5,0.5){$\xrightarrow[\textrm{\scriptsize 1 bot}]{\textrm{\scriptsize use}}$}
\rput(-4.5,0.5){$\xrightarrow[\textrm{\scriptsize 1 bot}]{\textrm{\scriptsize use}}$}
\rput(-0.5,0.5){$\xrightarrow[\textrm{\scriptsize 1 top}]{\textrm{\scriptsize use}}$}
\rput(4.5,0.5){$\xrightarrow[\textrm{\scriptsize 2 top}]{\textrm{\scriptsize use}}$}

\end{pspicture}
\end{center}

Now we describe the inverse of $\tau$. 
To prune a PM dot diagram means to delete the last dot in each row along with 
its incident edge and then join up their now-isolated partners 
unless the last dots were originally connected to each other in which 
case there is nothing to join up.
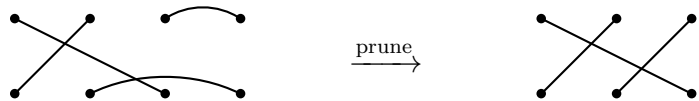
\begin{center} 
    
\begin{pspicture}(-5,0)(5,1)
    

\psbezier[linewidth=.8pt](-3,1)(-2.7,1.2)(-2.3,1.2)(-2,1)
\psbezier[linewidth=.8pt](-4,0)(-3.4,0.3)(-2.6,0.3)(-2,0)

\psdots(-5,1)(-5,0)(-4,1)(-4,0)(-3,1)(-3,0)(-2,1)(-2,0)
\psdots(4,1)(4,0)(3,1)(3,0)(2,1)(2,0)

\psline(-5,0)(-4,1)
\psline(-5,1)(-3,0)
\psline(2,1)(4,0)
\psline(2,0)(3,1)
\psline(3,0)(4,1)

\rput(0,.5){$\xrightarrow{\text{prune}}$ }

\end{pspicture}
\end{center} 

We will need the ``shift'' map $S$ that takes a dot that starts an 
upline in a PM dot diagram to one that does not:
if $i=i_{1}$ starts an upline then the 
upline ends at position $i_{2}>i_{1}$. If $i_{2}$ starts another 
upline, proceed similarly to get $i_{3}>i_{2}$ and so on, stopping at 
the first $i_{k}$ that does not start an upline, and set $S(i)=i_{k}$.
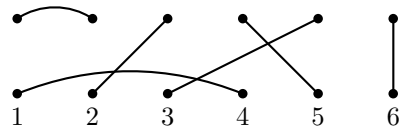
\begin{center} 
    \begin{pspicture}(-3,-2)(3,1.2)
    
\psbezier[linewidth=.8pt](-3,1)(-2.7,1.2)(-2.3,1.2)(-2,1)
\psbezier[linewidth=.8pt](-3,0)(-2,0.4)(-1,0.4)(0,0)

\psdots(-3,1)(-3,0)(-2,1)(-2,0)(-1,1)(-1,0)
\psdots(0,0)(0,1)(1,0)(1,1)(2,0)(2,1)

\psline(-2,0)(-1,1)
\psline(-1,0)(1,1)
\psline(0,1)(1,0)
\psline(2,0)(2,1)

 \rput(-3,-0.3){\textrm{{\footnotesize $1$}}}
\rput(-2,-0.3){\textrm{{\footnotesize $2$}}}
\rput(-1,-0.3){\textrm{{\footnotesize $3$}}}
\rput(0,-0.3){\textrm{{\footnotesize $4$}}}
\rput(1,-0.3){\textrm{{\footnotesize $5$}}}
\rput(2,-0.3){\textrm{{\footnotesize $6$}}}

\rput(0,-1.1){\small{There is an upline from 2 to 3 and another
from 3 to 5} }
\rput(0,-1.6){\small{but 5 does not start an upline. So $S(2)=5.$} }

\end{pspicture}
\end{center} 

At the step where a PM dot diagram of size $k$ is pruned to one of size $k-1$ record
the pair $(Y_{k},i_{k})$ of the build-matching code according as uplines
are created and/or destroyed in the pruning process. First, if the
last dots in each row are joined to each other, then
$(Y_{k},i_{k})=(B,k)$. Otherwise, let $i$ (resp. $j$) denote the
position in its row of the partner of the last dot in the top (resp.
bottom) row and consider four cases as in the Table.
\begin{center}
\begin{tabular}{|c|cccccc|}
    \hline
     & \raisebox{-0.3ex}[0pt]{row} & \raisebox{-0.3ex}[0pt]{row} &  & \raisebox{-0.3ex}[0pt]{upline}  &
     \raisebox{-0.3ex}[0pt]{upline} &   \\
    \raisebox{1.5ex}[0pt]{case}   & \raisebox{.3ex}[0pt]{of $i$} &
    \raisebox{.3ex}[0pt]{of $j$} &
    \raisebox{1.5ex}[0pt]{restriction} & \raisebox{.3ex}[0pt]{created?} &  \raisebox{.3ex}[0pt]{destroyed?} &
    \raisebox{1.5ex}[0pt]{\quad $(Y_{k},i_{k})$\quad} \\ \hline
    1a & top & top & - &   &   &   \\
    1b & top & bot & $i\le j$ & \raisebox{1.5ex}[0pt]{no} & \raisebox{1.5ex}[0pt]{no} &
    \raisebox{1.5ex}[0pt]{$(B,S(i))$}  \\  \hline
    {\rule[-3mm]{0mm}{8mm} 2}& top & bot & $i>j$ & yes & no & $(T,j)$  \\ \hline
    3a & bot & bot & - &  &  &   \\
    3b & bot & top & $i\ge j$ & \raisebox{1.5ex}[0pt]{no}  &
    \raisebox{1.5ex}[0pt]{yes}  & \raisebox{1.5ex}[0pt]{$(T,i)$}   \\ \hline
    {\rule[-3mm]{0mm}{8mm} 4} & bot & top & $i<j$ & yes & yes & $(B,i)$
    \\ \hline
\end{tabular}
\end{center}
An example of each case is shown.
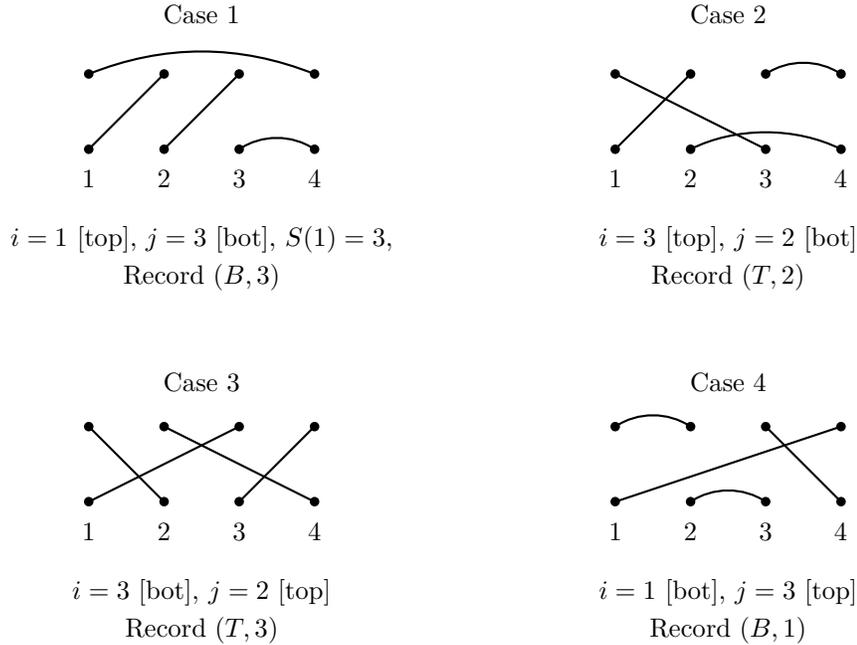
\begin{center}
\begin{pspicture}(-5,-1.7)(5,2)

\psbezier[linewidth=.8pt](-3,0)(-2.7,0.2)(-2.3,0.2)(-2,0)
\psbezier[linewidth=.8pt](-5,1)(-4,1.4)(-3,1.4)(-2,1)

\rput(-5,-0.4){\textrm{{\footnotesize $1$}}}
\rput(-4,-0.4){\textrm{{\footnotesize $2$}}}
\rput(-3,-0.4){\textrm{{\footnotesize $3$}}}
\rput(-2,-0.4){\textrm{{\footnotesize $4$}}}

\rput(2,-0.4){\textrm{{\footnotesize $1$}}}
\rput(3,-0.4){\textrm{{\footnotesize $2$}}}
\rput(4,-0.4){\textrm{{\footnotesize $3$}}}
\rput(5,-0.4){\textrm{{\footnotesize $4$}}}

\psdots(-3,1)(-3,0)(-2,1)(-2,0)(-4,1)(-4,0)(-5,1)(-5,0)
\psdots(3,1)(3,0)(2,1)(2,0)(4,1)(4,0)(5,1)(5,0)

\psline(-5,0)(-4,1)
\psline(-4,0)(-3,1)

\psbezier[linewidth=.8pt](3,0)(3.6,0.3)(4.4,0.3)(5,0)
\psbezier[linewidth=.8pt](4,1)(4.3,1.2)(4.7,1.2)(5,1)

\psline(2,0)(3,1)
\psline(2,1)(4,0)

\rput(-3.5,1.8){\textrm{{\footnotesize Case 1}}}
\rput(-3.5,-1.2){\textrm{{\footnotesize $i=1$ [top], $j=3$ [bot], $S(1)=3$,}}}
\rput(-3.5,-1.7){\textrm{{\footnotesize Record $(B,3)$}}}

\rput(3.5,1.8){\textrm{{\footnotesize Case 2}}}
\rput(3.5,-1.2){\textrm{{\footnotesize $i=3$ [top], $j=2$ [bot]}}}
\rput(3.5,-1.7){\textrm{{\footnotesize Record $(T,2)$}}}

\end{pspicture}
\end{center}

\begin{center}
\begin{pspicture}(-5,-1.7)(5,2.2)

\rput(-5,-0.4){\textrm{{\footnotesize $1$}}}
\rput(-4,-0.4){\textrm{{\footnotesize $2$}}}
\rput(-3,-0.4){\textrm{{\footnotesize $3$}}}
\rput(-2,-0.4){\textrm{{\footnotesize $4$}}}

\rput(2,-0.4){\textrm{{\footnotesize $1$}}}
\rput(3,-0.4){\textrm{{\footnotesize $2$}}}
\rput(4,-0.4){\textrm{{\footnotesize $3$}}}
\rput(5,-0.4){\textrm{{\footnotesize $4$}}}

\psdots(-3,1)(-3,0)(-2,1)(-2,0)(-4,1)(-4,0)(-5,1)(-5,0)
\psdots(3,1)(3,0)(2,1)(2,0)(4,1)(4,0)(5,1)(5,0)

\psline(-5,0)(-3,1)
\psline(-5,1)(-4,0)
\psline(-2,0)(-4,1)
\psline(-2,1)(-3,0)

\psbezier[linewidth=.8pt](2,1)(2.3,1.2)(2.7,1.2)(3,1)
\psbezier[linewidth=.8pt](3,0)(3.3,0.2)(3.7,0.2)(4,0)

\psline(2,0)(5,1)
\psline(4,1)(5,0)

\rput(-3.5,1.6){\textrm{{\footnotesize Case 3}}}
\rput(-3.5,-1.2){\textrm{{\footnotesize $i=3$ [bot], $j=2$ [top]}}}
\rput(-3.5,-1.7){\textrm{{\footnotesize  Record $(T,3)$}}}

\rput(3.5,1.6){\textrm{{\footnotesize Case 4}}}
\rput(3.5,-1.2){\textrm{{\footnotesize $i=1$ [bot], $j=3$ [top]}}}
\rput(3.5,-1.7){\textrm{{\footnotesize Record $(B,1)$}}}

\end{pspicture}
\end{center}
The entire pruning process for the \pm
$1\:3\,/\,2\:10\,/\,4\:7\,/\,5\:9\,/\,6\:8$ is shown:
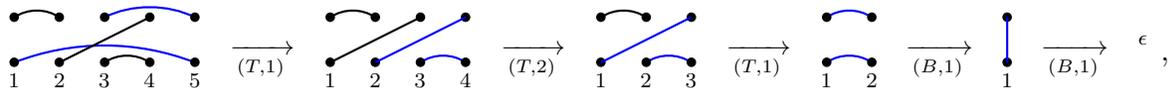
\begin{center}

\begin{pspicture}(-7,-0.2)(15,0.3)

\psset{xunit=.6cm,yunit=.6cm}

\psbezier[linecolor=blue,linewidth=.8pt](-11,0)(-9.7,0.5)(-8.3,0.5)(-7,0)
\psbezier[linewidth=.8pt](-9,0)(-8.7,0.2)(-8.3,0.2)(-8,0)
\psbezier[linewidth=.8pt](-11,1)(-10.7,1.2)(-10.3,1.2)(-10,1)
\psbezier[linecolor=blue,linewidth=.8pt](-9,1)(-8.4,1.3)(-7.6,1.3)(-7,1)
\psbezier[linewidth=.8pt](-4,1)(-3.7,1.2)(-3.3,1.2)(-3,1)
\psbezier[linecolor=blue,linewidth=.8pt](-2,0)(-1.7,0.2)(-1.3,0.2)(-1,0)

\psbezier[linecolor=blue,linewidth=.8pt](3,0)(3.3,0.2)(3.7,0.2)(4,0)
\psbezier[linecolor=blue,linewidth=.8pt](7,0)(7.3,0.2)(7.7,0.2)(8,0)
\psbezier[linewidth=.8pt](2,1)(2.3,1.2)(2.7,1.2)(3,1)
\psbezier[linecolor=blue,linewidth=.8pt](7,1)(7.3,1.2)(7.7,1.2)(8,1)

\psdots(-11,1)(-10,1)(-9,1)(-8,1)(-7,1)(-4,1)(-3,1)(-2,1)(-1,1)(2,1)(3,1)(4,1)(7,1)(8,1)(11,1)
\psdots(-11,0)(-10,0)(-9,0)(-8,0)(-7,0)(-4,0)(-3,0)(-2,0)(-1,0)(2,0)(3,0)(4,0)(7,0)(8,0)(11,0)

\psline(-10,0)(-8,1)
\psline(-4,0)(-2,1)
\psline[linecolor=blue](-3,0)(-1,1)
\psline[linecolor=blue](2,0)(4,1)
\psline[linecolor=blue](11,0)(11,1)

\rput(-5.5,0){\textrm{{\footnotesize  $\xrightarrow[(T,1)]{}$}}}
\rput(0.5,0){\textrm{{\footnotesize  $\xrightarrow[(T,2)]{}$}}}
\rput(5.5,0){\textrm{{\footnotesize  $\xrightarrow[(T,1)]{}$}}}
\rput(9.5,0){\textrm{{\footnotesize  $\xrightarrow[(B,1)]{}$}}}
\rput(12.5,0){\textrm{{\footnotesize  $\xrightarrow[(B,1)]{}$}}}

\rput(-11,-.4){\textrm{{\scriptsize 1}}}
\rput(-10,-.4){\textrm{{\scriptsize 2}}}
\rput(-9,-.4){\textrm{{\scriptsize 3}}}
\rput(-8,-.4){\textrm{{\scriptsize 4}}}
\rput(-7,-.4){\textrm{{\scriptsize 5}}}

\rput(-4,-.4){\textrm{{\scriptsize 1}}}
\rput(-3,-.4){\textrm{{\scriptsize 2}}}
\rput(-2,-.4){\textrm{{\scriptsize 3}}}
\rput(-1,-.4){\textrm{{\scriptsize 4}}}

\rput(2,-.4){\textrm{{\scriptsize 1}}}
\rput(3,-.4){\textrm{{\scriptsize 2}}}
\rput(4,-.4){\textrm{{\scriptsize 3}}}

\rput(7,-.4){\textrm{{\scriptsize 1}}}
\rput(8,-.4){\textrm{{\scriptsize 2}}}

\rput(11,-.4){\textrm{{\scriptsize 1}}}
\rput(14,.5){\textrm{{\scriptsize $\epsilon$}}}

\rput(14.5,0){,}

\end{pspicture}
\end{center}
yielding the build-matching code $\big( (B,1),(B,1),(T,1),(T,2),(T,1)\big)$.
\begin{prop}
    \label{tau}
    Under $\tau$, the pair $i\nearrow j$ is an upline iff $(T,i)$ occurs an odd 
    number of times in the build-matching 
    code and $j$ is the position in the code of the last occurrence 
    of a $(Y,i)\ (Y=T$ or $B)$.
\end{prop}
\noindent \emph{Proof}\quad The first occurrence of $T_{i}$ introduces an upline from $i$ that 
remains undisturbed by each later $Y_{j}$ with $j\ne i$, is switched 
to an upline from $i$ to the new top end dot by $B_{i}$, and is killed 
by a second occurrence of $T_{i}$. A third occurrence of $T_{i}$ 
reintroduces an upline from $i$. Thus an upline from $i$ is present in 
the resulting PM dot diagram iff $T_{i}$ occurs an odd number of times in the code 
and in this case the upline is $i\nearrow j$ where $j$ is the position 
of the last $Y_{i}$ in the code ($Y$ may be $B$ or $T$). \qed

That $\Phi$ sends Klazar violator/partner pairs $(i,j)$ to uplines $i\nearrow 
j$ follows from Props. \ref{sig} and \ref{tau}.

\section{Trapezoidal Words}

Riordan \cite[p.\,9]{rior76} considered the Cartesian product $[1]\times 
[3]\times [5]\times \ldots \times [2n-1]$. He called its entries, 
$(a_{k})_{1\le k \le n}$, \emph{trapezoidal words} and observed that 
the statistic ``number of distinct entries'' on trapezoidal words is 
distributed as what are now called the second-order Eulerian numbers,
\htmladdnormallink{A008517}{http://www.research.att.com:80/cgi-bin/access.cgi/as/njas/sequences/eisA.cgi?Anum=A008517}.
There is an obvious bijection from build-tree codes to trapezoidal 
words: $\big( (X_{k},i_{k})\big)_{1\le k \le n} \mapsto (a_{k})_{1\le k 
\le n}$ with $a_{k}=2i_{k}$ if $X_{k}=L$ and $a_{k}=2i_{k}+1$ if 
$X_{k}=R$. Props. \ref{sig} and \ref{tau} now yield the following corollary.
\begin{cor}
    \label{tri}
The following three statistics are equidistributed and all have the 
\gf 
\[
\left( \frac{1 - y}{2e^{x(y - 1)} - 1 - y } \right)^{\frac{1}{2}}
\]
of Cor. \ref{kvgf}.
\end{cor}

\vspace*{-5mm}

\begin{enumerate}
    \item ``\# Klazar violators'' on \os
    
    \item ``\# uplines'' on PM dot diagrams
    
    \item ``\# even entries that occur an odd number of times'' on 
    trapezoidal words.
\end{enumerate}

A variation of the bijection $\tau$ from build-matching codes to PM dot diagrams 
yields a generalization of the equidistribution of items 2 and 3 in 
Cor. \ref{tri}. Given a build-matching code, 
this time build up a PM dot diagram by successively enlarging the current PM dot diagram as 
follows. First, suppose $Y_{k}=B$. If there is a weak downline from $i$ use 
$i$ top, otherwise use the partner of $i$ (which may be top or bottom).
Now suppose $Y_{k}=T$. If there is an upline $i \nearrow j$ from $i$ use 
$j$ top, otherwise use $i$ bot. For example, the code 
$B_{1},T_{1},B_{2},T_{1}$ successively yields 
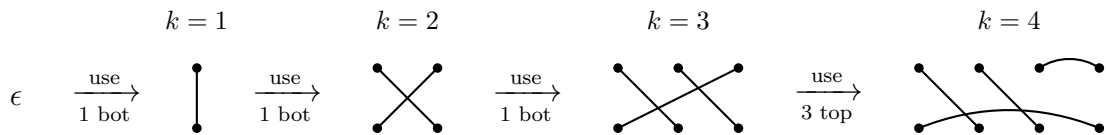
\begin{center} 
\begin{pspicture}(-8,0)(10,1.3)
\psset{unit=.8cm}
\psline(-6,0)(-6,1)
\psline(-3,0)(-2,1)
\psline(-2,0)(-3,1)
\psline(1,0)(3,1)
\psline(1,1)(2,0)
\psline(2,1)(3,0)
\psline(6,1)(7,0)
\psline(7,1)(8,0)

\psbezier[linewidth=.8pt](6,0)(7,0.4)(8,0.4)(9,0)
\psbezier[linewidth=.8pt](8,1)(8.3,1.2)(8.7,1.2)(9,1)

\psdots(-6,0)(-3,0)(-2,0)(1,0)(2,0)(3,0)(6,0)(7,0)(8,0)(9,0)
\psdots(-6,1)(-3,1)(-2,1)(1,1)(2,1)(3,1)(6,1)(7,1)(8,1)(9,1)

\rput(-6,1.8){\textrm{{\footnotesize $k=1$}}}
\rput(-2.5,1.8){\textrm{{\footnotesize $k=2$}}}
\rput(2,1.8){\textrm{{\footnotesize $k=3$}}}
\rput(7.5,1.8){\textrm{{\footnotesize $k=4$}}}

\rput(-9,0.5){$\epsilon$}

\rput(-7.5,0.5){$\xrightarrow[\textrm{\scriptsize 1 bot}]{\textrm{\scriptsize use}}$}
\rput(-4.5,0.5){$\xrightarrow[\textrm{\scriptsize 1 bot}]{\textrm{\scriptsize use}}$}
\rput(-0.5,0.5){$\xrightarrow[\textrm{\scriptsize 1 bot}]{\textrm{\scriptsize use}}$}
\rput(4.5,0.5){$\xrightarrow[\textrm{\scriptsize 3 top}]{\textrm{\scriptsize use}}$}

\end{pspicture}
\end{center}

\begin{cor}
    The joint distribution of the statistics ``\# even-to-odd matches'' and ``\# 
    odd-to-even matches'' on \ps of size $n$ is the same as that of the statistics 
    ``\# even entries that occur an odd number of times'' and 
    ``\# odd entries that occur an odd number of times'' on trapezoidal 
    words of length $n$.
\end{cor}
The \gf for ``\# uplines'' in PM dot diagrams is given in Cor. \ref{tri} 
and it is not hard to find the analogous \gf for ``\# vertical 
lines'': 
\[
\frac{1}{e^{x (1-y)} \sqrt{1-2x}}.
\]
Is there a nice \gf for the joint distribution of the three 
statistics ``\# uplines'', ``\# downlines'', and ``\# vertical 
lines'' in PM dot diagrams?

\end{document}